\pdfoutput=1
\documentclass[9pt,twoside,epsf]{article}
\usepackage{stmaryrd}
\usepackage{amsmath,latexsym,amssymb,amsfonts,amsbsy}
\usepackage{bm}
\usepackage{graphicx,epstopdf}
\usepackage[caption=false]{subfig}
\usepackage{cases}
\usepackage{multirow}
\usepackage{booktabs}
\usepackage{mathdots}

\newfont{\bb}{msbm10}
\def\Bbb#1{\mbox{\bb #1}}

\def\T{\top}

\def\diag{{\rm diag}}

\def\rank{{\rm rank}}

\def\VEC{{\rm vec}}

\usepackage{colortbl}
\usepackage{ragged2e}
\usepackage[margin=2em,labelsep=space,skip=0.5em,font=normalfont]{caption}
\usepackage{algorithm}  
\usepackage{algpseudocode} 
\floatname{algorithm}{\color{black} Algorithm}
\algrenewcommand{\algorithmiccomment}[1]{\quad{\color{red}\%\ #1}}
\numberwithin{algorithm}{section}


\newtheorem{method}{Method}[section]

\newtheorem{remark}{Remark}[section]
\newtheorem{theorem}{Theorem}[section]

\newtheorem{lemma}{Lemma}[section]

\newcommand{\reals}{\makebox{{\Bbb R}}}
\newcommand{\ceals}{\makebox{{\Bbb C}}}

\newcommand{\Z}{\makebox{{\Bbb Z}}}

\newcommand{\bbb}[1]{\text{\bf #1}}

\makeatletter 

\@addtoreset{equation}{section}
\makeatother 

\begin{document}
\cleardoublepage
\pagestyle{myheadings}
\bibliographystyle{plain}

\title{Sine-transform-based fast solvers for Riesz fractional nonlinear Schrödinger equations with attractive nonlinearities
\thanks{Corresponding author: Xi Yang (yangxi@nuaa.edu.cn); Other authors: Chao Chen (chenchao3@nuaa.edu.cn), Fei-Yan Zhang (zhangfy@nuaa.edu.cn).}}

\author{
Chao Chen\thanks{School of Mathematics, Nanjing University of Aeronautics and Astronautics, Nanjing 211106, China.}\ \thanks{Key Laboratory of Mathematical Modelling and High Performance Computing of Air Vehicles (NUAA), MIIT, Nanjing 211106, China.},
Xi Yang\footnotemark[2]\ \footnotemark[3],
Fei-Yan Zhang\thanks{Wuxi Taihu University, Wuxi 214064, China.}
}

\maketitle

\begin{abstract}
This paper presents fast solvers for linear systems arising from the discretization of fractional nonlinear Schrödinger equations with Riesz derivatives and attractive nonlinearities. These systems are characterized by complex symmetry, indefiniteness, and a $d$-level Toeplitz-plus-diagonal structure. We propose a Toeplitz-based anti-symmetric and normal splitting iteration method for the equivalent real block linear systems, ensuring unconditional convergence. The derived optimal parameter is approximately equal to 1. By combining this iteration method with sine-transform-based preconditioning, we introduce a novel preconditioner that enhances the convergence rate of Krylov subspace methods. Both theoretical and numerical analyses demonstrate that the new preconditioner exhibits a parameter-free property (allowing the iteration parameter to be fixed at 1). The eigenvalues of the preconditioned system matrix are nearly clustered in a small neighborhood around 1, and the convergence rate of the corresponding preconditioned GMRES method is independent of the spatial mesh size and the fractional order of the Riesz derivatives.

	\bigskip
	
	
	\bigskip
	
	
\end{abstract}

\section{Introduction}\label{intro}
The Schrödinger equation is one of the most important equations in quantum mechanics. As is well known, the standard Schrödinger equation (SSE) is a second-order partial differential equation that can be derived from the path integral over the Brownian motion. The fractional Schrödinger equation (FSE) can be obtained by extending the standard diffusion operator in SSE to the fractional diffusion operator. There are two ways to do this: one is to generalize to the fractional Laplacian \cite{Laskin01,Laskin02}, and the other is to generalize to the Riesz fractional derivative \cite{1997chaos}. Given the crucial role of  FSE in quantum mechanics, its theoretical properties and practical applications have been extensively studied \cite{2012BSS,2013Luchko,2006Guo,2008HanGuo,2015Duo,2019Li}. However, the nonlocal nature of fractional derivatives poses challenges in obtaining exact solutions for FSE. Numerical methods have become crucial in studying FSE. In recent decades, several classes of numerical methods have been established for FSE, such as finite element methods \cite{LM2018JCP,LM2017NUMA}, spectral methods \cite{DSW2016CMA,WY2019ANM}, collocation methods \cite{Amore2010JMP,Bhrawy2017ANM}, and finite difference methods \cite{WDL2013JCP,WDL2014JCP,WPD2015JCP,ZRP2019SCM,2022dis2D}, etc.

In this paper, we consider the following Riesz fractional nonlinear Schrödinger equation (RFNSE)
\begin{equation} \label{rfnse}
    {\imath}u_{t}+\frac{\partial^{\alpha}u(\bbb{x},t)}{\partial|\bbb{x}|^{\alpha}}+\rho|u|^2 u=0, \quad
	 \bbb{x} \in\reals^{d}, 0 < t \le \mathbf{t}
\end{equation}
with the initial condition
	$$	u(\bbb{x},0)=u_{0}(\bbb{x}),\quad \bbb{x} \in \reals^{d} , $$
where $\imath=\sqrt{-1}$, $\bbb{x}=(x_{1},x_{2},\dots,x_{d})^{\T} \in \reals^{d}$, $\rho>0$ is a real constant (representing attractive interaction of particles), and $u_{0}(\bbb{x})$ is a complex-valued function. The Riesz fractional derivative $\partial^{\alpha}u(\bbb{x},t)/\partial|\bbb{x}|^{\alpha}$ with $1<\alpha\leq 2$ defined by
\begin{align*}
	\frac{\partial^{\alpha}u(\bbb{x},t)}{\partial|\bbb{x}|^{\alpha}}=
-\frac{1}{2\text{cos}(\frac{\alpha \pi}{2})}\sum_{j=1}^{d}[_{-\infty}D_{x_{j}}^{\alpha}u(\bbb{x},t)+_{x_{j}}D_{\infty}^{\alpha}u(\bbb{x},t)],
\end{align*}
in which the left and right Riemann–Liouville fractional derivatives are defined by
\begin{align*}
	_{-\infty}D_{x_{j}}^{\alpha}u(\bbb{x},t)&=\frac{1}{\Gamma(2-\alpha)}\frac{\partial^{2}}{\partial x_{j}^{2}}\int_{-\infty}^{x_{j}}\frac{u(x_{1}, \dots, u_{j-1}, \xi, u_{j+1}, \dots, x_{d},t)}{(x_{j}-\xi)^{\alpha-1}}\text{d}\xi,\\
	_{x_{j}}D_{\infty}^{\alpha}u(\bbb{x},t)&=\frac{1}{\Gamma(2-\alpha)}\frac{\partial^{2}}{\partial x_{j}^{2}}\int_{x_{j}}^{\infty}\frac{u(x_{1}, \dots, u_{j-1}, \xi, u_{j+1}, \dots, x_{d},t)}{(\xi-x_{j})^{\alpha-1}}\text{d}\xi,
\end{align*}
where $\Gamma(\cdot)$ is the gamma function. When $\alpha=2$, the above RFNSE reduces to the standard nonlinear Schrödinger equation (SNSE)  \cite{Pan2017CMA,Pan2020ANM}.

If a linearly implicit difference scheme is applied to discretize RFNSE (\ref{rfnse}), a sequence of complex linear systems of the form $(D_{d}-T_{d}+\imath I)\bbb{u}=\bbb{b}$ is obtained, where $D_{d}$ is a diagonal matrix, and $T_{d}$ is a symmetric positive definite $d$-level Toeplitz matrix. Moreover, the parameter $\rho>0$ causes $D_{d}$ to be a positive semi-definite matrix, thus making $D_{d}-T_{d}$ an indefinite matrix. Therefore, the indefinite coefficient matrix $D_{d}-T_{d}+\imath I$ can be considered as a complex symmetric matrix or a $d$-level Toeplitz-plus-diagonal matrix.

Since the above complex linear systems have a $d$-level Toeplitz-plus-diagonal structure, fast direct solvers are not available for them. However, for $d$-level Toeplitz matrices, their matrix-vector multiplication can be achieved through the fast Fourier transform (FFT) \cite{FFT}. Hence, the Krylov subspace iteration methods can be reasonable options to solve the complex linear systems derived from (\ref{rfnse}). It is worth noting that an ill-conditioned coefficient matrix will result in high computational cost and slow convergence rate. In recent years, preconditioning techniques for ($d$-level) Toeplitz-plus-diagonal linear systems have been extensively studied to improve the convergence rate of Krylov subspace iteration methods. For instance, Chan and Ng introduced an effective banded preconditioner to solve Toeplitz-plus-band linear systems (allowing the Toeplitz-plus-diagonal structure as a special case) \cite{ChanNg1993}. Ng and Pan proposed an approximate inverse circulant-plus-diagonal (AICD) preconditioner for Toeplitz-plus-diagonal matrices \cite{NgPan2010}. Bai et al. studied the diagonal and Toeplitz splitting (DTS) preconditioner for 1-level \cite{BaiLuPan2017} and $d$-level \cite{BaiLu2020} Toeplitz-plus-diagonal systems derived from one-dimensional (1D) and higher-dimensional fractional diffusion equations, significantly improving the convergence rate of Krylov subspace iteration methods. Furthermore, Lu et al. combined the DTS preconditioner with the sine-transform to propose a new efficient preconditioner \cite{Lu2021JAMC}, further reducing the computational cost of Krylov subspace iteration methods for Toeplitz-plus-diagonal linear systems. However, the aforementioned preconditioners are primarily designed for Hermitian positive definite ($d$-level) Toeplitz-plus-diagonal systems, rather than complex symmetric and indefinite ones.

The matrix $D_{d}-T_{d}+\imath I$ also has a complex symmetric structure, which can be handled by three classes of methods. The first is the class of alternating-type iteration methods. Bai et al. proposed the Hermitian and skew-Hermitian splitting (HSS \cite{HSS10}) iteration method and the preconditioned HSS (PHSS \cite{PHSS2004}) iteration method for non-Hermitian positive definite matrices. For complex symmetric matrices with positive semi-definite real and imaginary parts, and at least one of them being positive definite, Bai et al. introduced the modified HSS (MHSS \cite{MHSS18}) iteration method and the preconditioned MHSS (PMHSS \cite{PMHSS2011,PMHSS2013}) iteration method. The second is the class of preconditioned Krylov subspace iteration methods, whose convergence rates can be significantly improved by carefully designed preconditioners. For instance, preconditioners can be naturally derived from the alternating-type iteration methods. The third is the class of C-to-R iteration methods \cite{Axelsson2000}, which can be interpreted as applying the block Gaussian elimination to the real block 2-by-2 linear system equivalent to the complex symmetric one. The related Schur-complement linear subsystem requires high-quality and efficient preconditioners, which are often difficult to construct. Unfortunately, the aforementioned methods may not be suitable or may perform poorly for complex symmetric indefinite linear systems.

In this paper, we consider an equivalent real block 2-by-2 form of $(D_{d}-T_{d}+\imath I)\bbb{u}=\bbb{b}$, and construct the Toeplitz-based anti-symmetric and normal (TBAN) splitting. Combining the spirit of the alternating direction implicit iteration \cite{ADI1955,ADI}, we propose the TBAN iteration method. Theoretically, we prove that the TBAN iteration method converges unconditionally for any positive iteration parameter, and the optimal iteration parameters is deducted. This new iteration method naturally leads to the TBAN preconditioner. The implementation of this preconditioner requires to solve two linear subsystems with coefficient matrices $ \bigl[ \begin{smallmatrix} \omega I & T_{d} \\	-T_{d} & \omega I \end{smallmatrix} \bigr]$ and $ \bigl[\begin{smallmatrix} (\omega+1) I & -D_{d} \\ D_{d} & (\omega+1) I \end{smallmatrix} \bigr]$. It is well known that circulant preconditioning works well for linear systems with Toeplitz structures \cite{Lei2016AEJAM,Lei2013JCP}. However, the behavior of Krylov subspace iteration methods with circulant-preconditioning is linearly dependent on the system size \cite{CC1992}. In fact, when the Toeplitz matrix is real symmetric, it can also be well approximated by the $\tau$-matrix \cite{Benedetto1995,Benedetto1990,1983taudiagonal,2022spec,Huang2021ANM}. Compared to the circulant-preconditioning, the $\tau$-preconditioning (i.e., the sine-transform-based preconditioning) reduces the computational cost by a constant factor at each iteration, and has a better convergence rate \cite{Benedetto1990}. Therefore, we combine the TBAN preconditioner with the $\tau$-preconditioning to construct a new preconditioner, which can be efficiently implemented by the fast sine-transform (FST) \cite{1983taudiagonal,2023CMA}. Theoretically and numerically, we have proved that the new sine-transform-based preconditioner has parameter-free property, with the optimal parameter approximately equal to 1, and exhibiting better performance than the circulant-based preconditioner. The eigenvalues of the preconditioned system matrix are clustered in a small neighborhood around 1, and the convergence rate of the preconditioned GMRES method is independent of the spatial mesh size and the fractional order of Riesz derivatives.

This paper is organized as follows. In Section \ref{discre}, the complex linear systems $(D_{d}-T_{d}+\imath I)\bbb{u}=\bbb{b}$ are derived by applying linearly implicit difference schemes to RFNSE, see \cite{WDL2014JCP} for the 1D case and \cite{2022dis2D} for the two-dimensional (2D) case. In Section \ref{iter}, we introduce the TBAN iteration method, prove its unconditional convergence property, and deduct the optimal iteration parameter. In Section \ref{precon}, we propose a novel preconditioner based on the sine-transform, and analyze the eigenvalue distribution of the preconditioned system matrix. Numerical results are presented in Section \ref{exp} to illustrate the reliability and efficiency of the new preconditioner. Finally, concluding remarks are given in Section \ref{clu}.

\section{The discrete linear systems}\label{discre}
In practical computations, the $d$-dimensional space $\reals^{d}$ is truncated to a bounded domain $\Omega$, and the problem (\ref{rfnse}) is equipped with the Dirichlet boundary condition
$$u(\bbb{x},t)=0, \quad \bbb{x}\in
\partial\Omega, \quad  0 \le t \le \mathbf{t}.$$
The discretization of (\ref{rfnse}) using linearly implicit difference schemes results in the aforementioned complex linear system \begin{align} \label{equ3}
	(D_{d}-T_{d}+\imath I)\ \bbb{u} &= \bbb{b}
\end{align}
at each time level. Although the new method presented in this paper can be applied to linear systems derived from the fractional Schrödinger equation (\ref{rfnse}) for $d=1,2,3$, we only consider the 1D (\cite{WDL2014JCP}) and 2D (\cite{2022dis2D}) cases. This is because, to the best of our knowledge, the work on stable and conservative linearly implicit difference schemes for the three-dimensional (3D) fractional Schrödinger equation (\ref{rfnse}) is currently unavailable in the literature.

\subsection{The 1D case}\label{sec-1D-DLS}
Let $d=1$ and $\Omega=[\mbox{a},\mbox{b}]$. Given positive integers $N$ and $M$, the space interval and time interval can be divided into $M+1$ and $N$ equal parts respectively, i.e., the spatial step size $h=(\text{b}-\text{a})/(M+1)$, and the temporal step size $\Delta t=\mathbf{t}/N$. Let $u_{j}^{n} \approx u(x_{j},t_{n})$ represent the numerical solution at the spatial position $x_{j}=a+jh$ and the time level $t_{n}=n\Delta t$.

The 1D Riesz fractional derivative can be discretized by the fractional centered difference scheme \cite{2012JCP,2006Riesz} in the bounded interval $[\mbox{a},\mbox{b}]$ as
\begin{gather*}
	\frac{\partial^{\alpha}}{\partial|x|^{\alpha}}u(x_{j},t)
	=-\frac{1}{h^{\alpha}}\sum^M_{k=1}c_{j-k}u_{k}+\mathcal{O}(h^2),
\end{gather*}
where the coefficients $c_{k}$ read as
\begin{gather} \label{coefficients}
	c_{k}=\frac{(-1)^{k}\Gamma(\alpha+1)}{\Gamma(\alpha/2-k+1) \Gamma(\alpha/2+k+1)}, \forall\ k\in\Z,
\end{gather}
and they satisfy \cite{2012JCP}
\begin{gather} \nonumber
	c_{0}\ge 0,\
	c_{k}=c_{-k}\le0\ (\forall\ k\ge 1),\ \mbox{and}\
	\sum^{+\infty}_{k=-\infty,k\ne0}|c_{k}|=c_{0}.
\end{gather}

The linearly implicit conservative difference (LICD) scheme \cite{WDL2014JCP} applied to the 1D truncated problem (\ref{rfnse}) results in
\begin{equation}  \label{discretizedRFSE}
		\imath\frac{u^{n+1}_j-u^{n-1}_j}{2\Delta t}-\frac{1}{h^\alpha} \sum^M_{k=1}c_{j-k}\hat{u}^n_k+\rho|u^n_j|^2 \hat{u}^n_j=0,
\end{equation}
where $\hat{u}^n_j=(u^{n+1}_j+u^{n-1}_j)/2$, for $j=1,2,\dots,M,n=1,2,\dots,N-1$, and the initial and boundary conditions read as
 $u^0_j=u_0(x_j)$, $u^{n}_0=u^{n}_{M+1}=0$. The matrix-vector form of (\ref{discretizedRFSE}) is
\begin{align} \label{discretizedCNLSMaxForm}
	(D_{1}^{n+1}-T_{1}+\imath I)\bbb{u}^{n+1} &= \bbb{b}^{n+1}, \forall\ n\ge1.
\end{align}
Here, $\bbb{u}^{n+1}=[u_1^{n+1},\ldots,u_M^{n+1}]^{\T}$, $D_{1}^{n+1}=\diag\{d_1^{n+1},d_2^{n+1},\ldots,d_M^{n+1}\}\in\reals^{M\times M}$ ($d_j^{n+1}=\rho\Delta t|u_j^{n}|^2$) is a diagonal matrix, $I\in\reals^{M\times M}$ is an identity matrix, and $T_{1}=\mu T_0$ is the symmetric Toeplitz matrix with $\mu=\Delta t / h^\alpha$ and
\begin{align}\label{equ2}
	T_0 &=
	\begin{bmatrix}
		c_0 & c_{-1} & \ldots & c_{2-M} & c_{1-M} \\
		c_1 & c_0 & \ddots & \ddots & c_{2-M} \\
		\vdots & \ddots & \ddots & \ddots & \vdots \\
		c_{M-2} & \ddots & \ddots & c_0 & c_{-1} \\
		c_{M-1} & c_{M-2} & \ldots & c_1 & c_0 \\
	\end{bmatrix}.
\end{align}

\subsection{The 2D case}
Let $d=2$ and $\Omega =[\mbox{a},\mbox{b}]\times[\mbox{a}',\mbox{b}']$. For positive integers $N$ and $M$, let $\Delta t=\mathbf{t}/N$, $h_x=(\mbox{b}-\mbox{a})/(M+1)$, $h_y=(\mbox{b}'-\mbox{a}')/(M+1)$ be the temporal step size and the spatial step sizes. Then, the space-time domain $\Omega \times [0,\mathbf{t}]$ can be covered by
$$\Omega_{h}=\left\{(x_{j},y_{k},t_{n})|x_{j}=\mbox{a}+jh_x, y_{k}=\mbox{a}'+kh_y, t_{n}=n\Delta t, j,k=0,1,\ldots M+1, n=0,\ldots ,N\right\}.
$$
Additionally, let $u_{j,k}^{n} \approx u(x_{j},y_{k},t_{n})$ be the numerical
solution at the spatial position $(x_{j},y_{k})$ and the time level $t_{n}$.

The 2D Riesz fractional derivative can be approximated by the fractional centered difference scheme as follows
\begin{gather*}
	\left(\frac{\partial^{\alpha}}{\partial|x|^{\alpha}}+\frac{\partial^{\alpha}}{\partial|y|^{\alpha}}\right)u(x_{j}, y_k,t)
	=-\frac{1}{h_x^{\alpha}}\sum^M_{s=1}c_{j-s}u_{s,k}-\frac{1}{h_y^{\alpha}}\sum^M_{s=1}c_{k-s}u_{j,s}+\mathcal{O}(h^2),
\end{gather*}
where $c_{k}$ is defined in (\ref{coefficients}).

The three-level linearized implicit difference (TLID) scheme \cite{2022dis2D} applied to the 2D truncated problem (\ref{rfnse}) results in
\begin{equation}  \label{discretizedRFSE2D}
	\imath\frac{u_{j,k}^{n+1}-u^{n-1}_{j,k}}{2\Delta t}-\frac{1}{h_x^\alpha} \sum^M_{s=1}c_{j-s}\hat{u}^n_{s,k}-\frac{1}{h_y^\alpha} \sum^M_{s=1}c_{k-s}\hat{u}^n_{j,s}+\rho|u^n_{j,k}|^2 \hat{u}^n_{j,k}=0,
\end{equation}
where $\hat{u}^n_{j,k}=(u_{j,k}^{n+1}+u^{n-1}_{j,k})/2$, for $j,k=1,2,\dots,M,n=1,2,\dots,N-1$, and the initial and boundary conditions are
\begin{align*}
  u_{j,k}^{0}&=u_{0}(x_{j},y_{k}),\quad j,k=1,\ldots,M, \\
  u_{0,k}^{n}&=u_{M+1,k}^{n}=u_{j,0}^{n}=u_{j,M+1}^{n}=0,\quad n=0,\ldots,N.
\end{align*}
The matrix-vector form of (\ref{discretizedRFSE2D}) is
\begin{align} \label{discretizedRFSE2DMaxForm}
	(D_{2}^{n+1}-T_{2}+\imath I)\bbb{u}^{n+1} &= \bbb{b}^{n+1},\ \forall\ n\ge 1,
\end{align}
where $\bbb{u}^{n+1}=\VEC(U^{n+1})\in \ceals^{M^{2}}$ is the unknown vector (`$\VEC$' stacks the columns of a matrix, and $U^{n+1}=[u_{j,k}^{n+1}] \in \ceals^{M\times M}$), $D_{2}^{n+1} = \rho\Delta t\, \diag\{\VEC(\bar{U}^n \odot U^n)\} \in \reals^{M^{2} \times M^{2}}$ is a diagonal matrix (`$\odot$' represents the Hadamard product, and $U^{n}=[u_{j,k}^{n}] \in \ceals^{M\times M}$), $I\in \reals^{M^{2} \times M^{2}}$ is an identity matrix, and $T_2$ is the 2-level Toeplitz matrix of the form
\begin{align}\label{2DToeplitz}
	T_{2}=I\otimes T_{x} +T_{y}\otimes I\in\reals^{M^{2} \times M^{2}}
\end{align}
with $T_{x}=\mu_x T_0$ ($\mu_x=\Delta t / h_x^\alpha$) and $T_{y}=\mu_y T_0$ ($\mu_y=\Delta t / h_y^\alpha$).

\section{The TBAN iteration method}\label{iter}

Considering the linear system (\ref{equ3}), let the solution be $\bbb{u}=y+\imath z $, and the right-hand side be $\bbb{b}=p+\imath q$, where $y$, $z$, $p$ and $q$ are real vectors. Then, (\ref{equ3}) can be equivalently rewritten as the following real non-symmetric positive definite block linear system
\begin{align}\label{positiveBlockForm}
	\mathcal{R}_{d}x &\equiv
	\begin{bmatrix}
		I & T_{d}-D_{d}\\
		D_{d}-T_{d} & I\\
	\end{bmatrix}
	\begin{bmatrix}
		z \\
		y\\
	\end{bmatrix}
	=
	\begin{bmatrix}
		-p \\
		q \\
	\end{bmatrix}
	\equiv f.
\end{align}
By disrupting the $d$-level Toeplitz-plus-diagonal block $T_{d}-D_{d}$ in (\ref{positiveBlockForm}), the system matrix $\mathcal{R}_{d}$ admits a Toeplitz-based anti-symmetric and normal  (TBAN) splitting of the form
\begin{align} \label{TBANsplitting}
	\mathcal{R}_{d} &= \begin{bmatrix}
		0 & T_{d}\\
		-T_{d} & 0\\
	\end{bmatrix} + \begin{bmatrix}
	I & -D_{d}\\
	D_{d} & I\\
	\end{bmatrix} \equiv \mathcal{T}_{d} + \mathcal{D}_{d},
\end{align}
where $\mathcal{T}_{d}$ is an anti-symmetric matrix with $d$-level Toeplitz block $T_{d}$, and $\mathcal{D}_{d}$ is a normal matrix. By combining the TBAN splitting (\ref{TBANsplitting}) and the spirit of the alternating direction implicit iteration \cite{ADI1955,ADI}, we propose the following TBAN iteration method.
\begin{method}[The TBAN iteration method]
	\label{NATBiterationMd}
	Let $x^{(0)} $ be an arbitrary initial guess. For $k=0,1,2,\ldots$ until the sequence of iterates $\{x^{(k)}\}_{k\ge 0}$ converges, compute the next iterate $x^{(k+1)}$ according to the following procedure:
	\begin{align}\label{NATBiterationSc}
		\left\{\begin{aligned}
			(\omega I+\mathcal{T}_{d})x^{(k+\frac{1}{2})} & = (\omega I-\mathcal{D}_{d})x^{(k)}+f, \\
			(\omega I+\mathcal{D}_{d})x^{(k+1)} & = (\omega I-\mathcal{T}_{d})x^{(k+\frac{1}{2})}+f,
		\end{aligned}\right.
	\end{align}
	where $\omega>0$ is a prescribed positive constant.
\end{method}

	
Through simple calculations, the TBAN iteration (\ref{NATBiterationSc}) can be integrated into the following fixed-point iteration scheme
$$\mathcal{F}_{d,\omega}x^{(k+1)} = \mathcal{G}_{d,\omega}x^{(k)}+f,$$
where the matrices
\begin{align}\label{defpreconditionger1}
	\left\{\begin{aligned}
		\mathcal{F}_{d,\omega} = \frac{1}{2\omega}
		(\omega I+\mathcal{T}_{d})(\omega I+\mathcal{D}_{d}), \\ \mathcal{G}_{d,\omega} = \frac{1}{2\omega}
		(\omega I-\mathcal{T}_{d})(\omega I-\mathcal{D}_{d}),
	\end{aligned}\right.
\end{align}
constitute the following splitting
   \begin{align}\label{NATBsplitting2}
	\mathcal{R}_{d} = \mathcal{F}_{d,\omega}-\mathcal{G}_{d,\omega}.
\end{align}
The iteration matrix of the TBAN iteration (\ref{NATBiterationSc}) is
\begin{align} \label{NATBiterMatrix}
	\mathcal{L}_{d,\omega} &= \mathcal{F}_{d,\omega}^{-1}\mathcal{G}_{d,\omega}.
\end{align}
	
The following theorem gives the convergence property of the TBAN iteration method.
	
\begin{theorem} \label{NATBconvergenceThm}
Let $\mathcal{R}_{d}$ be a non-symmetric positive definite block matrix as defined in (\ref{positiveBlockForm}). Let $\mathcal{T}_{d}$ and $\mathcal{D}_{d}$ constitute a TBAN splitting of $\mathcal{R}_{d}$ in (\ref{TBANsplitting}). Let $\omega$ be a positive constant. For any initial vector $x^{(0)}$, the TBAN iteration sequence $\{x^{(k)}\}_{k\ge 0}$ converges to the unique solution of the linear system (\ref{positiveBlockForm}). Furthermore, the spectral radius of the TBAN iteration matrix $\rho(\mathcal{L}_{d,\omega})$ is bounded as follows
		\begin{align} \label{NASSrateBoundedByOne}
			\rho(\mathcal{L}_{d,\omega}) \le \sigma(\omega) < 1,\ \forall\ \omega > 0,
		\end{align}
		with
		\begin{align} \label{NATBupperBound}
			\sigma(\omega) &= \sqrt{\frac {(\omega-1)^2+\lambda_{\max}^2}{(\omega+1)^2+\lambda_{\max}^2}},
		\end{align}
		where $\lambda_{\max}$ is the maximum diagonal of $D_{d}$.
	\end{theorem}
	
{\em Proof.}
Obviously, for $\omega>0$, both the matrices $\omega I+\mathcal{T}_{d}$ and $\omega I+\mathcal{D}_{d}$ are positive definite. We notice the following relation
	$$\widehat{\mathcal{L}}_{d,\omega}
	\equiv (\omega I+ \mathcal{D}_{d})\mathcal{L}_{d,\omega}(\omega I +\mathcal{D}_{d})^{-1}=\mathcal{U}_{d,\omega} \mathcal{V}_{d,\omega}$$
	with $\mathcal{U}_{d,\omega} = (\omega I+\mathcal{T}_{d})^{-1}(\omega I-\mathcal{T}_{d})$ and $\mathcal{V}_{d,\omega} = (\omega I-\mathcal{D}_{d})(\omega I+\mathcal{D}_{d})^{-1}$, thus it holds that
	\begin{align} \nonumber
		\rho(\mathcal{L}_{d,\omega}) &= \rho(\mathcal{U}_{d,\omega} \mathcal{V}_{d,\omega}) \\ \label{NATBrateBoundedBy2norm}
		&\le \|\mathcal{U}_{d,\omega}\|_2 \|\mathcal{V}_{d,\omega}\|_2.
	\end{align}
	
Next, we will estimate the upper bounds of $\|\mathcal{U}_{d,\omega}\|_2$ and $\|\mathcal{V}_{d,\omega}\|_2$.
\begin{itemize}
  \item For the bound of $\|\mathcal{U}_{d,\omega}\|_2$, based on the fact that $\mathcal{T}^{\T}_{d}\mathcal{T}_{d}=\mathcal{T}_{d}\mathcal{T}_{d}^{\T}$, we have
	\begin{align*} \nonumber
		\mathcal{U}_{d,\omega}^{\T}\mathcal{U}_{d,\omega} &= (\omega I -\mathcal{T}_{d})^{\T}(\omega I +\mathcal{T}_{d})^{-\T}(\omega I +\mathcal{T}_{d})^{-1}(\omega I -\mathcal{T}_{d}) \\ \nonumber
		&= (\omega I +\mathcal{T}_{d})^{-\T}(\omega I +\mathcal{T}_{d})^{-1}(\omega I -\mathcal{T}_{d})^{\T}(\omega I -\mathcal{T}_{d}) \\ \nonumber
		&= [(\omega I +\mathcal{T}_{d})(\omega I +\mathcal{T}_{d})^{\T}]^{-1}(\omega I -\mathcal{T}_{d})^{\T}(\omega I -\mathcal{T}_{d})\\
		&= \begin{bmatrix}
			\omega^2 I +T_{d}^2& 0 \\
			0 & \omega^2 I +T_{d}^2 \
		\end{bmatrix}^{-1}
		\begin{bmatrix}
			\omega^2 I +T_{d}^2 & 0 \\
			0 & \omega^2 I +T_{d}^2
		\end{bmatrix} \\
		&= I.
	\end{align*}
	Therefore, $\mathcal{U}_{d,\omega}$ is an orthogonal matrix, i.e.,
		\begin{align} \label{NATBupBoundFactorU}
		\left \|\mathcal{U}_{d,\omega}\right \|_2 &=   1.
	\end{align}
  \item For the bound of $\|\mathcal{V}_{d,\omega}\|_2$, the fact $\mathcal{D}_{d}^{\T}\mathcal{D}_{d}=\mathcal{D}_{d}\mathcal{D}_{d}^{\T}$ leads to
	\begin{align*}
	\mathcal{V}_{d,\omega}^{\T}\mathcal{V}_{d,\omega} &= (\omega I+\mathcal{D}_{d}^{\T})^{-1} (\omega I-\mathcal{D}_{d}^{\T}) (\omega I-\mathcal{D}_{d}) (\omega I+\mathcal{D}_{d})^{-1} \\
		&= (\omega I-\mathcal{D}_{d}^{\T}) (\omega I-\mathcal{D}_{d}) [(\omega I+\mathcal{D}_{d}) (\omega I+\mathcal{D}_{d}^{\T})]^{-1} \\
		&=
		\begin{bmatrix}
			(\omega-1)^2 I+D_{d}^2 & 0 \\
			0 & (\omega-1)^2 I+D_{d}^2 \
		\end{bmatrix}
		\begin{bmatrix}
			(\omega+1)^2 I+D_{d}^2 & 0 \\
			0 & (\omega+1)^2 I+D_{d}^2 \
		\end{bmatrix}^{-1}.\\
	\end{align*}
	Since $D_{d}$ is a positive semi-definite diagonal matrix, it reads that
	\begin{align}\label{norm-2-v}
		\|\mathcal{V}_{d,\omega}\|_2 &= \max_{\lambda_i\in\lambda(D_{d})} g(\omega;\lambda_i)
                                        = g(\omega;\lambda_{\max}) < 1,
	\end{align}
    where $\lambda(D_{d})$ contains the diagonals of $D_{d}$, and $g(\omega;\lambda)=\sqrt{ [(\omega-1)^2+\lambda^2] / [(\omega+1)^2+\lambda^2]}$. The second equality of (\ref{norm-2-v}) is due to the increasing monotonicity of $g(\omega;\lambda)$ with respect to $\lambda>0$.
\end{itemize}
	
	By combining the relations (\ref{NATBrateBoundedBy2norm})-(\ref{norm-2-v}), we can conclude that $\rho(\mathcal{L}_{d,\omega}) \le \sigma(\omega) <1$ holds true for $\omega>0$.	

	$\hfill\square$

	\begin{remark}\label{remarkTBANiter}
We provide some remarks regarding the optimal parameter and the related optimal convergence rate of the TBAN iteration.
		\begin{enumerate}
			\item The optimal parameter $\omega^{\star}$ minimizing $\sigma(\omega)$ can be obtained by determining the positive root of the equation $\rm{d}[\sigma(\omega)] / \rm{d} \omega=0$, i.e.,
				$$\omega^{\star}=\sqrt{\lambda_{\max}^2+1}.$$
				By adopting the optimal parameter $\omega^{\star}$, the convergence rate of the TBAN iteration satisfies
				\begin{align*}
					\rho(\mathcal{L}_{d,\omega^{\star}}) \le \sigma(\omega^{\star}) =\frac{\lambda_{\max}}{1+\sqrt{\lambda_{\max}^2+1}}.
				\end{align*}
				\item If the solution $u(\mathbf{x},t)$ of the truncated fractional Schrödinger equation (\ref{rfnse}) is uniformly bounded, then for any given parameter $\rho>0$, it follows that $\lambda_{\max}=\mathcal{O}(\Delta t)$ due to the structure of the diagonals of $D_{d}$. Thus, we have
$$\omega^{\star}=\sqrt{\mathcal{O}(\Delta t)^2+1} \approx 1 \ \text{and} \ \rho(\mathcal{L}_{d,\omega^{\star}})\le \sigma(\omega^{\star}) \approx \mathcal{O}(\Delta t)$$
for a small temporal step size $\Delta t$. In addition, if we adopt the parameter $\omega=1$, according to (\ref{NATBupperBound}), the eigenvalues of the TBAN iteration matrix $\mathcal{L}_{d,1}$ stay in the neighborhood of $0$ with the radius
$$\frac{\lambda_{\max}}{\sqrt{4+\lambda_{\max}^2}}
=\frac{\mathcal{O}(\Delta t)}{\sqrt{4+\mathcal{O}(\Delta t)^2}}
= \mathcal{O}(\Delta t).$$
				\end{enumerate}
					
			\end{remark}

\section{Preconditioning}\label{precon}
For the sake of brevity, we will focus on the case where $h_x=h_y=h$. The analysis for the case where $h_x \neq h_y$ is not significantly more complicated.

The matrix $\mathcal{F}_{d,\omega}$ in the TBAN splitting (\ref{NATBsplitting2}) can naturally serve as a preconditioner, called the TBAN preconditioner, of the linear system (\ref{positiveBlockForm}). The main workload for implementing $\mathcal{F}_{d,\omega}$ in the preconditioned Krylov subspace iteration methods lies in solving the following generalized residual (GR) linear systems
\begin{align} \label{generalizedResEqns}
	\mathcal{F}_{d,\omega}z^{(k)} = r^{(k)},\ \forall\ k\ge 0,
\end{align}
where $r^{(k)}$ is the current residual vector at the $k$-th iteration, and $z^{(k)}$ is the corresponding GR vector. Specifically, two linear subsystems with coefficient matrices $\omega I+\mathcal{T}_{d}$  and $\omega I+\mathcal{D}_{d}$ need to be solved. To reduce the computational costs, an improved version of $\mathcal{F}_{d,\omega}$ based on the sine-transform is considered.

For the 1-level Toeplitz matrix $T_1$ with its first row being $\mu [c_{0}, c_{1}, \ldots , c_{M-1}]$, a natural way to construct the related $\tau$-matrix approximation $\tau(T_1)$ is given as below
\begin{align}\label{taudefine}
	\tau(T_1) = T_1 - \text{HC}(T_1),
\end{align}
where $\text{HC}(T_1)$ is the Hankel correction of $T_1$ \cite{1999taudefine}, i.e.,
\begin{align} \label{HCform}
	\text{HC}(T_1)=\mu \left[\begin{array}{cccccccc}
		c_2 & \cdots & c_{M-1} & 0 & 0 \\
		\vdots & \iddots & 0 & 0 & 0 &  \\
		c_{M-1} & \iddots & \iddots & \iddots & c_{M-1} \\
		0 & 0 & 0 & \iddots & \vdots  \\
		0 & 0 & c_{M-1} & \cdots & c_{2}
	\end{array}\right].
\end{align}
It is worth noting that $\tau(T_1)$ can be diagonalized by the sine-transform \cite{1983taudiagonal,1999taudefine}, i.e.,
\begin{align} \label{taudiagonal}
	\tau(T_1)=S \Lambda S,
\end{align}
where $\Lambda$ is diagonal holding all the eigenvalues of $\tau(T_1)$ determined by its first column, and $S$ is symmetric with the following elements
\begin{align}
	[S]_{i,j}=\sqrt{\frac{1}{M+1}}\text{sin}\left(\frac{\pi ij}{M+1}\right), \quad 1\le i,j \le M.
\end{align}
In addition, the matrix-vector multiplication of $S$ with a vector can be computed with $\mathcal{O}(M \log M)$ operations by FST \cite{1983taudiagonal,2023CMA}.

The sine-transform-based preconditioner can be constructed by approximating the $d$-level Toeplitz matrix in $\mathcal{F}_{d,\omega}$ in the following way
\begin{align}\label{taupre}
	\widetilde{\mathcal{F}}_{d,\omega} = \frac{1}{2\omega}
	\left[\omega I+\tau(\mathcal{T}_{d})\right](\omega I+\mathcal{D}_{d}), \quad d=1, 2.
\end{align}
Here, $d=1$ and $d=2$ represent the 1D case and the 2D case respectively, and
$\tau(\mathcal{T}_{d})= \bigl[ \begin{smallmatrix} 0 & \tau(T_{d}) \\
	-\tau(T_{d}) & 0 \end{smallmatrix} \bigr]$ with $\tau(T_{1})$ being defined in (\ref{taudefine}) and
	\begin{align}\label{2Dtau}
		\tau(T_{2}) = I \otimes \tau(T_{1}) + \tau(T_{1}) \otimes I.
	\end{align}
	
Now we focus on the 2D case. Specifically, we study the eigenvalue clustering of the 2D preconditioned system matrix $\widetilde{\mathcal{F}}^{-1}_{2,\omega}\mathcal{R}_{2}$. Since $\widetilde{\mathcal{F}}^{-1}_{2,\omega} \mathcal{R}_{2} $ can be factorized as
\begin{align} \label{NATBpreconSysMax}
	\widetilde{\mathcal{F}}^{-1}_{2,\omega} \mathcal{R}_{2} &= \underbrace{\widetilde{\mathcal{F}}^{-1}_{2,\omega} \mathcal{F}_{2,\omega}}\ \underbrace{\mathcal{F}^{-1}_{2,\omega} \mathcal{R}_{2}},
\end{align}
we need to study the eigenvalue distribution of $\widetilde{\mathcal{F}}^{-1}_{2,\omega} \mathcal{F}_{2,\omega}$ and $\mathcal{F}^{-1}_{2,\omega} \mathcal{R}_{2}$.

Firstly, the following theorem gives the eigenvalue distribution of $\mathcal{F}^{-1}_{2,\omega} \mathcal{R}_{2}$ .
\begin{theorem}\label{theorem1}
	Let $\mathcal{R}_{2} \in\reals^{2M^{2}\times 2M^{2}}$ be the system matrix in (\ref{TBANsplitting}), $\mathcal{F}_{2,\omega} $ be the 2D TBAN preconditioner, $\omega$ be a positive constant, and $\sigma(\omega)$ be defined by (\ref{NATBupperBound}). Then, the eigenvalues of the 2D TBAN preconditioned system matrix $\mathcal{F}_{2,\omega}^{-1}\mathcal{R}_{2}$ are located in a circle of radius $\sigma(\omega)<1$ centered at $1$.
\end{theorem}

{\em Proof.}
	It can be easily proved based on the fact $\mathcal{F}_{2,\omega}^{-1} \mathcal{R}_{2}= I - \mathcal{L}_{2,\omega}$ derived from (\ref{NATBsplitting2}) and (\ref{NATBiterMatrix}), and the result of Theorem \ref{NATBconvergenceThm}.
$\hfill\square$

\begin{remark}\label{remark41}
	Theorem \ref{theorem1} and Remark \ref{remarkTBANiter} indicate that when $\omega=\omega^{\star}$ or $\omega =1$, the eigenvalues of $\mathcal{F}_{2,\omega}^{-1} \mathcal{R}_{2}$ are situated within a circle of radius $\mathcal{O}(\Delta t)$ centered at $1$, particularly in the case of a small temporal step $\Delta t$.
\end{remark}

Secondly, we study the eigenvalue distribution of $\widetilde{\mathcal{F}}^{-1}_{2,\omega} \mathcal{F}_{2,\omega}$.
The following Lemmas \ref{lemma2}-\ref{lemma5} provide the eigenvalue bounds for $T_{1}$, $T_{2}$, $\tau(T_{1})$, $\tau(T_{2})$.

\begin{lemma}\label{lemma1}(\cite{2023Yang}) Let
	$$\theta=\frac{\left(1-\frac{1+\alpha}{5+\alpha / 2}\right)^{5+\frac{\alpha}{2}} e^{1+\alpha} \Gamma(\alpha+1) \sin \left(\frac{\pi \alpha}{2}\right)}{\pi \alpha} \quad \text { and } \quad \theta_0=\frac{\sqrt{2} e^{13 / 12} \Gamma(\alpha+1) \sin \left(\frac{\pi \alpha}{2}\right)}{\pi \alpha}$$ be two constants, $k_0 \geq 3$, and $1<\alpha \leq 2$. Then,
	$$
    \frac{\theta}{(k_{0}+1/2)^{\alpha}}<\sum_{j=k_0+1}^{\infty}\left|c_j\right|<\frac{\theta_{0}}{(k_{0}-1)^{\alpha}},	$$
    where  $c_j=(-1)^j \Gamma(\alpha+1) /[\Gamma(\alpha / 2-j+1) \Gamma(\alpha / 2+j+1)]$.
\end{lemma}

\begin{lemma}(\cite{2023Yang})\label{lemma2}
	Consider the 1-level Toeplitz matrix $T_{1}\in\reals^{M \times M}$, let $M$ be even, $1<\alpha \leq 2$. Then,
	$$
	\frac{2  \Delta t \theta}{(\mathrm{b}-\mathrm{a})^\alpha}<\lambda_{T_{1}}<\frac{2  \Delta t}{h^\alpha}\left[\frac{\Gamma(\alpha+1)}{\Gamma(\alpha / 2+1)^2}-\frac{\theta h^\alpha}{(\mathrm{b}-\mathrm{a})^\alpha}\right], \quad M \geq 4,
	$$
	where $\lambda_{T_{1}}$ represents any eigenvalue of $T_{1}$.
\end{lemma}

\begin{lemma}\label{lemma3}
	Let $M$ be even, $1<\alpha \leq 2$. The eigenvalues of the 2-level Toeplitz matrix $T_{2}=I\otimes T_{1} +T_{1}\otimes I\in\reals^{M^{2} \times M^{2}}$
	satisfy
	$$	\frac{4  \Delta t \theta}{(\mathrm{b}-\mathrm{a})^\alpha}< \lambda_{T_{2}} <\frac{4  \Delta t}{h^\alpha}\left[\frac{\Gamma(\alpha+1)}{\Gamma(\alpha / 2+1)^2}-\frac{\theta h^\alpha}{(\mathrm{b}-\mathrm{a})^\alpha}\right], \quad M \geq 4,$$
	where $\lambda_{T_{2}}$ represents any eigenvalue of $T_{2}$.
\end{lemma}

{\em Proof.}
The above bounds can be directly obtained from Lemma \ref{lemma2} and the property of the Kronecker product.
$\hfill\square$

\begin{lemma}\label{lemma4}
	Let $\tau(T_{1})$ be the $\tau$-matrix approximation of $T_{1}$ defined in (\ref{taudefine}). Let $M$ be even, $1<\alpha \leq 2$. Then, the eigenvalues of $\tau(T_{1})$ are bounded as
	\begin{align*}
		\frac{2\Delta t\theta}{(b-a)^{\alpha}}&< \lambda_{\tau(T_{1})} <\frac{2\Delta t}{h^{\alpha}}\left[\frac{\Gamma(\alpha+1)}{\Gamma(\alpha/2+1)^{2}}-\frac{\theta h^{\alpha}}{(b-a)^{\alpha}}\right]-\frac{c_{2}\Delta t}{h^{\alpha}}, \quad M \geq 4,
	\end{align*}
where $\lambda_{\tau(T_{1})}$ represents any eigenvalue of $\tau(T_{1})$.
\end{lemma}

{\em Proof.}
Let $R_{i}[\cdot]$ be the sum of absolute values of the off-diagonal entries in the $i$-th row of a matrix, and $a_{ii}\left[\cdot\right]$ denotes the $i$-th diagonal entry of a matrix. According to the decaying property of $c_k=(-1)^k \Gamma(\alpha+1) /[\Gamma(\alpha / 2-k+1) \Gamma(\alpha / 2+k+1)]$, it is easy to get

$$R_{1}[T_{1}] < R_{2}[T_{1}] < \cdots <R_{\frac{M}{2}}[T_{1}]=R_{\frac{M}{2}+1}[T_{1}],$$
and
$$R_{\frac{M}{2}}[T_{1}]=R_{\frac{M}{2}+1}[T_{1}] > R_{\frac{M}{2}+2}[T_{1}] > \cdots >R_{M}[T_{1}].$$
Since the off-diagonal elements of $T_{1}$ and HC$(T_{1})$ are all negative, we have
\begin{align*}
	\max_{i}R_{i}[\tau(T_{1})] < R_{\frac{M}{2}}[T_{1}] < 2\mu\sum_{k=1}^{M-1}\left | c_{k} \right |.
\end{align*}
Besides, we have
\begin{align*}
	\min_{i}a_{ii}\left[\tau(T_{1})\right]=\mu c_{0} \quad \text{and} \quad  \max_{i}a_{ii}\left[\tau(T_{1})\right]=\mu(c_{0}-c_{2}).
\end{align*}
According to Gerschgorin disk theorem and Lemma \ref{lemma1}, it reads that
\begin{align*}
	\mu(c_{0}-2\sum_{k=1}^{M-1}\lvert c_{k} \rvert)&\leq \lambda_{\tau(T_{1})} \leq\mu(c_{0}-c_{2}+2\sum_{k=1}^{M-1}\lvert c_{k} \rvert))
	\\	2\mu\sum_{k=M}^{+\infty}\lvert c_{k} \rvert &\leq \lambda_{\tau(T_{1})} \leq 2\mu(c_{0}-\sum_{k=M}^{+\infty}\lvert c_{k} \rvert)-\mu c_2
	\\   \frac{2\Delta t\theta}{(b-a)^{\alpha}}&< \lambda_{\tau(T_{1})} <\frac{2\Delta t}{h^{\alpha}}\left[\frac{\Gamma(\alpha+1)}{\Gamma(\alpha/2+1)^{2}}-\frac{\theta h^{\alpha}}{(b-a)^{\alpha}}\right]-\frac{c_{2}\Delta t}{h^{\alpha}}, \quad  M>4
\end{align*}
$\hfill\square$
\begin{lemma}\label{lemma5}
	Let $\tau(T_{2})$ be the $\tau$-matrix approximation of $T_{2}$ defined in (\ref{2Dtau}). Let $M$ be even, $1<\alpha \leq 2$. Then, the eigenvalues of $\tau(T_{2})$ are bounded as
	$$	\frac{4\Delta t\theta}{(b-a)^{\alpha}}<\lambda(\tau(T_{2}))<\frac{4\Delta t}{h^{\alpha}}\left[\frac{\Gamma(\alpha+1)}{\Gamma(\alpha/2+1)^{2}}-\frac{\theta h^{\alpha}}{(b-a)^{\alpha}}\right]-\frac{2c_{2}\Delta t}{h^{\alpha}}, \quad M \geq 4, $$
	where $\lambda_{\tau(T_{2})}$ represents any eigenvalue of $\tau(T_{2})$.
\end{lemma}

{\em Proof.}
The above bounds can be directly obtained from Lemma \ref{lemma4} and the property of the Kronecker product.
$\hfill\square$

Lemmas \ref{lemma3} and \ref{lemma5} show that $T_{2}$ and $\tau(T_{2})$ are positive definite. The following Lemma \ref{lemma6} depicts the extent to which $\tau(T_{2})$ approximates $T_{2}$.

\begin{lemma}\label{lemma6}
	Let $T_{2}\in\reals^{M^{2} \times M^{2}}$  and  $\tau(T_{2}) \in\reals^{M^{2} \times M^{2}}$ be defined in (\ref{2DToeplitz}) and(\ref{2Dtau}), respectively. Let $\epsilon$ be a small positive constant satisfying $2^{2\alpha +1} \mu \theta_0 / M^{\alpha}<\epsilon \le 2\mu \theta_{0}$, $M>4$ be even, and $k_0 = \lceil (2\mu \theta_{0} / \epsilon)^{1/\alpha} \rceil +1 $, where $\lceil$.$\rceil $ represents rounding a real number to positive infinity. Then, there exist two matrices, $U\in \mathbb{R}^{ M^{2} \times  M^{2}}$ and $V\in \mathbb{R}^{ M^{2} \times  M^{2}}$, satisfying that
	\begin{align*}
		T_{2} -\tau(T_{2})= U+V,
	\end{align*}
	where $\operatorname{rank} (U)<4M(k_{0}-1)$,
	$$\|U\|_{\infty}<2\mu \left[\frac{c_0}{2}-\frac{\theta}{\left(M-\frac{1}{2}\right)^\alpha}\right]\quad \text{and} \quad \|V\|_{\infty}< \epsilon.$$
\end{lemma}

{\em Proof.}
We split the Hankel correction HC$(T_{1})$ as
\begin{gather*}
	\text{HC}(T_{1})= \hat{U}+\hat{V},\ \text{with}\ \hat{U} = \mu\begin{bmatrix}
		c_{2} & \cdots & c_{k_{0}} & 0 & \cdots & 0 \\
		\vdots & \iddots & \iddots &  &  & \vdots \\
		c_{k_{0}} & \iddots &  &   &  & 0 \\
        0 & & & & \iddots & c_{k_{0}} \\
		\vdots &  &  & \iddots & \iddots & \vdots \\
		0 & \cdots & 0 & c_{k_{0}} & \cdots & c_{2}
	\end{bmatrix}.
\end{gather*}
Obviously, it holds that $\rank(\hat{U}) = 2(k_{0}-1)$ with $2\le k_{0} < 1+M/4$. Thanks to Lemma \ref{lemma1}, the $\ell_{\infty}$-norm estimates of $\hat{U}$ and $\hat{V}$ reads that
$$
\begin{aligned}
	\|\hat{U}\|_{\infty}
	& \leq \mu \sum_{j=1}^{M-1}\left|c_j\right|=\mu\left(\frac{c_0}{2}-\sum_{j=M}^{\infty}\left|c_j\right|\right) \\
	&< \mu \left[\frac{c_0}{2}-\frac{\theta}{\left(M-\frac{1}{2}\right)^\alpha}\right],
\end{aligned}
$$

$$
\begin{aligned}
	\|\hat{V}\|_{\infty} &  = \mu \sum_{j=k_0+1}^{M-1} \left |c_j\right |<\mu \sum_{j=k_0+1}^{\infty}\left|c_j\right| \\
	& <\frac{\mu \theta_0}{\left(k_0-1\right)^\alpha}<\frac{\epsilon}{2}.
\end{aligned}
$$
Straight forward computations lead to
\begin{align*}
	T_{2}-\tau(T_{2})&=I\otimes T_{1}+T_{1}\otimes I- I\otimes \tau(T_{1})-\tau(T_{1})\otimes I \\
	&=I\otimes (T_{1}-\tau(T_{1}))+(T_{1}-\tau(T_{1})) \otimes I \\
	&=I\otimes(\hat{U}+\hat{V})+(\hat{U}+\hat{V})\otimes I \\
    &=U+V,
\end{align*}
where $U=I\otimes \hat{U}+\hat{U}\otimes I$ and $V=I\otimes \hat{V}+\hat{V}\otimes I$. Then, $\rank(U)$ is bounded as follows
\begin{align*}
	\operatorname{rank} (U)&\le \operatorname{rank} (I\otimes \hat{U})+\operatorname{rank} (\hat{U}\otimes I) \\
	&=2\operatorname{rank}(\hat{U})\operatorname{rank}(I) \\
	&<4M(k_{0}-1),
\end{align*}
the $\ell_{\infty}$-norm estimate of $U$ satisfies
\begin{align*}
	\|U\|_{\infty}&\le\|I\otimes \hat{U}\|_{\infty}+\|\hat{U}\otimes I\|_{\infty}
	=2\|\hat{U}\|_{\infty}\|I\|_{\infty}
	<2\mu \left[\frac{c_0}{2}-\frac{\theta}{\left(M-\frac{1}{2}\right)^\alpha}\right],
\end{align*}
and the $\ell_{\infty}$-norm estimate of $V$ can be achieved as below
\begin{align*}
	\|V\|_{\infty}\le 2\|\hat{V}\|_{\infty}\|I\|_{\infty} <\epsilon.
\end{align*}
$\hfill\square$

The following Theorem \ref{theorem2} provides the property of $\widetilde{\mathcal{F}}^{-1}_{2,\omega} \mathcal{F}_{2,\omega}$.
\begin{theorem}\label{theorem2}
	Let $\mathcal{F}_{2,\omega}\in \mathbb{R}^{2 M^{2} \times 2 M^{2}}$ and $\widetilde{\mathcal{F}}_{2,\omega}\in \mathbb{R}^{2 M^{2} \times 2 M^{2}}$ be defined in (\ref{defpreconditionger1}) and (\ref{taupre}), respectively. Let $\epsilon$ be a small positive constant satisfying $2^{2\alpha +1} \mu \theta_0 / M^{\alpha}<\epsilon \le 2\mu \theta_{0}$, $M>4$ be even, $k_0 = \lceil (2\mu \theta_{0} / \epsilon)^{1/\alpha} \rceil +1 $, and $\nu = \max_{\mu_i \in \lambda(D_{2})} |\mu_i|$. Then, there exist two matrices, $U_{\tau}\in \mathbb{R}^{2 M^{2} \times 2 M^{2}}$ and $V_{\tau}\in \mathbb{R}^{2 M^{2} \times 2 M^{2}}$, satisfying that
	\begin{align} \label{tildeFinvF}
		\widetilde{\mathcal{F}}^{-1}_{2,\omega} \mathcal{F}_{2,\omega} &= I + U_{\tau}+V_{\tau},
	\end{align}
	where $\operatorname{rank} (U_{\tau})<8M(k_{0}-1)$,
	$$\left\| U_{\tau}\right\|_2<\frac{2\left[(\omega+1)^{2}+\nu^{2}\right]^{\frac{1}{2}}M^{\frac{1}{2}}\mu}{(\omega+1)\sqrt{\omega^{2}+\left[\frac{4\Delta t\theta}{(b-a)^{\alpha}}\right]^{2}}}\left[\frac{c_0}{2}-\frac{\theta}{\left(M-\frac{1}{2}\right)^\alpha}\right] \quad \text{and} \quad \left\| V_{\tau}\right\|_2<\frac{\left[(\omega+1)^{2}+\nu^{2}\right]^{\frac{1}{2}}M^{\frac{1}{2}}\epsilon}{(\omega+1)\sqrt{\omega^{2}+\left[\frac{4\Delta t\theta}{(b-a)^{\alpha}}\right]^{2}}}.$$
\end{theorem}

{\em Proof.}
Simple calculations lead to
\begin{align*}
	\widetilde{\mathcal{F}}_{2,\omega}^{-1}\mathcal{F}_{2,\omega}-I
	&=(\omega I+\mathcal{D}_{2})^{-1}[\omega I+\tau(\mathcal{T}_{2})]^{-1}(\omega I+\mathcal{T}_{2})(\omega I+\mathcal{D}_{2})-I \\
	&=(\omega I+\mathcal{D}_{2})^{-1}\left\{[\omega I+\tau(\mathcal{T}_{2})]^{-1}(\omega I+\mathcal{T}_{2})-I\right\} (\omega I+\mathcal{D}_{2}) \\
	&=(\omega I+\mathcal{D}_{2})^{-1} [\omega I+\tau(\mathcal{T}_{2})]^{-1} [\mathcal{T}_{2}-\tau(\mathcal{T}_{2}) ](\omega I+\mathcal{D}_{2}) \\
	&=(\omega I+\mathcal{D}_{2})^{-1}(\omega I+\tau(\mathcal{T}_{2}))^{-1}\begin{bmatrix}
		0 & U+V \\
		-(U+V) & 0
	\end{bmatrix}(\omega I+\mathcal{D}_{2}) \\
	&=U_{\tau}+V_{\tau},
\end{align*}
where$$U_{\tau}=(\omega I+\mathcal{D}_{2})^{-1}(\omega I+\tau(\mathcal{T}_{2}))^{-1}\begin{bmatrix}
	0 & U \\
	-U & 0
\end{bmatrix}(\omega I+\mathcal{D}_{2})$$ and $$V_{\tau}=(\omega I+\mathcal{D}_{2})^{-1}(\omega I+\tau(\mathcal{T}_{2}))^{-1}\begin{bmatrix}
	0 & V \\
	-V & 0
\end{bmatrix}(\omega I+\mathcal{D}_{2}).$$
Thus, according to Lemmas \ref{lemma5} and \ref{lemma6}, we have $\operatorname{rank}(U_{\tau})<8M(k_{0}-1)$, the $\ell_2$-norm estimate of $U_{\tau}$ reading that
\begin{align*}
	\left\| U_{\tau}\right\|_2 &\leq \left\| \begin{bmatrix}
		(\omega+1) I & -D_{2} \\
		D_{2} & (\omega+1) I
	\end{bmatrix}^{-1} \right\|_{2} \left\| \begin{bmatrix}
		(\omega+1) I & -D_{2} \\
		D_{2} & (\omega+1) I
	\end{bmatrix}\right\|_{2} \left\| \begin{bmatrix}
		\omega I & \tau(T_{2}) \\
		-\tau(T_{2}) & \omega I
	\end{bmatrix}^{-1} \right\|_{2} \left\| U \right\|_{2} \\
	&=\frac{\max\limits_{\mu_{i}\in\lambda(D_{2})}\sqrt{(\omega+1)^{2}+\mu_{i}^{2}}}{\min\limits_{\mu_{i}\in\lambda(D_{2})}\sqrt{(\omega+1)^{2}+\mu_{i}^{2}}\mathop{\min}\limits_{\lambda_{i}\in\lambda(\tau(T_{2}))}\sqrt{\omega^{2}+\lambda_{i}^{2}}}\left\| U \right\|_{2} \\
	&< \frac{\max\limits_{\mu_{i}\in\lambda(D_{2})}\sqrt{(\omega+1)^{2}+\mu_{i}^{2}}}{\min\limits_{\mu_{i}\in\lambda(D_{2})}\sqrt{(\omega+1)^{2}+\mu_{i}^{2}}\min\limits_{\lambda_{i}\in\lambda(\tau(T_{2}))}\sqrt{\omega^{2}+\lambda_{i}^{2}}}M^{\frac{1}{2}} \| U \|_\infty \\
	&<\frac{2\left[(\omega+1)^{2}+\nu^{2}\right]^{\frac{1}{2}}M^{\frac{1}{2}}\mu}{(\omega+1)\sqrt{\omega^{2}+\left[\frac{4\Delta t\theta}{(b-a)^{\alpha}}\right]^{2}}}\left[\frac{c_0}{2}-\frac{\theta}{\left(M-\frac{1}{2}\right)^\alpha}\right],
\end{align*}
and the $\ell_2$-norm estimate of $V_{\tau}$ reading that
\begin{align*}
	\|V_{\tau}\|_{2} <\frac{\left[(\omega+1)^{2}+\nu^{2}\right]^{\frac{1}{2}}M^{\frac{1}{2}}\epsilon}{(\omega+1)\sqrt{\omega^{2}+\left[\frac{4\Delta t\theta}{(b-a)^{\alpha}}\right]^{2}}}.
\end{align*}
$\hfill\square$

Finally, we investigate the property of $\widetilde{\mathcal{F}}^{-1}_{2,\omega}\mathcal{R}_{2}$ as follows.

\begin{theorem} \label{theorem3}
	Let $\mathcal{R}_{2}\in \mathbb{R}^{2 M^{2} \times 2 M^{2}}$ be defined in (\ref{positiveBlockForm}). Let $\epsilon$ be a small positive constant satisfying $2^{2\alpha +1} \mu \theta_0 / M^{\alpha}<\epsilon \le 2\mu \theta_{0}$, $M>4$ be even, $k_0 = \lceil (2\mu \theta_{0} / \epsilon)^{1/\alpha} \rceil +1 $, and $\nu = \max_{\mu_i \in \lambda(D_{2})} |\mu_i|$. Then, there exist two matrices, $\mathcal{P}_\tau \in \mathbb{R}^{2 M^{2} \times 2 M^{2}}$ and $\mathcal{Q}_\tau \in \mathbb{R}^{2 M^{2} \times 2 M^{2}}$, satisfying that
$$
\widetilde{\mathcal{F}}^{-1}_{2,\omega} \mathcal{R}_{2}=\mathcal{F}_{2,\omega}^{-1} \mathcal{R}_{2}+\mathcal{P}_\tau+\mathcal{Q}_\tau,
$$
where $\operatorname{rank} (\mathcal{P}_{\tau})<8M(k_{0}-1)$,
$$\|\mathcal{P}_{\tau}\|_{2}<\frac{4\left[(\omega+1)^{2}+\nu^{2}\right]M^{\frac{1}{2}}\mu}{(\omega+1)^{2}\sqrt{\omega^{2}+\left[\frac{4\Delta t\theta}{(b-a)^{\alpha}}\right]^{2}}}\left[\frac{c_0}{2}-\frac{\theta}{\left(M-\frac{1}{2}\right)^\alpha}\right] \quad \text{and} \quad \|\mathcal{Q}_{\tau}\|_{2}<\frac{2\left[(\omega+1)^2 + \nu^2\right]M^{\frac{1}{2}} \epsilon}{(\omega+1)^2\sqrt{\omega^2+[\frac{4 \Delta t \theta}{(\rm{b}-\rm{a})^\alpha}]^2}}.$$
\end{theorem}

{\em Proof.}
Due to the byproducts $\left\|\mathcal{U}_{2,\omega}\right\|_2=1$, $\left\|\mathcal{V}_{2,\omega}\right\|_2<1$ in the proof of Theorem \ref{NATBconvergenceThm}, and the fact
\begin{align} \nonumber
	\mathcal{F}_{2,\omega}^{-1} \mathcal{R}_{2}
	&\nonumber=(\omega I+\mathcal{D}_{2})^{-1}\left(I-\mathcal{U}_{2,\omega}\mathcal{V}_{2,\omega}\right)(\omega I+\mathcal{D}_{2}),
\end{align}
it reads that $\left\|I-\mathcal{U}_{2,\omega}\mathcal{V}_{2,\omega}\right\|_2 \leq\|I\|_2+\left\|\mathcal{U}_{2,\omega}\right\|_2\left\|\mathcal{V}_{2,\omega}\right\|_2<2$.
Then, the $\ell_2$-norm estimate of $\mathcal{F}_{2,\omega}^{-1} \mathcal{R}_{2}$ holds that
\begin{align*}
\|\mathcal{F}_{2,\omega}^{-1} \mathcal{R}_{2}\|_{2}\le \frac{2\sqrt{(\omega+1)^2 + \nu^2}}{\omega+1}.
\end{align*}
From (\ref{NATBpreconSysMax}) and (\ref{tildeFinvF}), we know that
\begin{align*}
	\widetilde{\mathcal{F}}^{-1}_{2,\omega} \mathcal{R}_{2} & =\left(I+U_{\tau}+V_{\tau}\right) \mathcal{F}_{2,\omega}^{-1} \mathcal{R}_{2} \\
	& =\mathcal{F}_{2,\omega}^{-1} \mathcal{R}_{2}+\mathcal{P}_\tau+\mathcal{Q}_\tau,
\end{align*}
where $\mathcal{P}_\tau=U_{\tau} \mathcal{F}_{2,\omega}^{-1} \mathcal{R}_{2}$, and $\mathcal{Q}_\tau=V_{\tau} \mathcal{F}_{2,\omega}^{-1} \mathcal{R}_{2}$. According to Theorem \ref{theorem2}, we have $\operatorname{rank}\left(\mathcal{P}_\tau\right)<8M (k_0-1)$,
\begin{align}
	\left\|\mathcal{P}_\tau\right\|_2 & \nonumber \leq\left\|U_{\tau}\right\|_2\left\|\mathcal{F}_{2,\omega}^{-1} \mathcal{R}_{2}\right\|_2 \\
	&  \nonumber <   \frac{4\left[(\omega+1)^{2}+\nu^{2}\right]M^{\frac{1}{2}}\mu}{(\omega+1)^{2}\sqrt{\omega^{2}+\left[\frac{4\Delta t\theta}{(b-a)^{\alpha}}\right]^{2}}}\left[\frac{c_0}{2}-\frac{\theta}{\left(M-\frac{1}{2}\right)^\alpha}\right],
\intertext{and}\nonumber
	\left\|\mathcal{Q}_\tau\right\|_2 & \leq\left\|V_{\tau}\right\|_2\left\|\mathcal{F}_{2,\omega}^{-1} \mathcal{R}_{2}\right\|_2 \\
	&  \nonumber <  \frac{2\left[(\omega+1)^2 + \nu^2\right]M^{\frac{1}{2}} \epsilon}{(\omega+1)^2\sqrt{\omega^2+[\frac{4 \Delta t \theta}{(\rm{b}-\rm{a})^\alpha}]^2}}.
\end{align}
$\hfill\square$

\begin{remark}
According to Theorem \ref{theorem3}, $\mathcal{Q}_\tau$ has a small norm, which implies that the eigenvalues of $\mathcal{F}_{2,\omega}^{-1} \mathcal{R}_{2}+\mathcal{Q}_\tau$ remain within the union of the $\mathcal{O}(\epsilon)$-neighborhoods of the eigenvalues of  $\mathcal{F}_{2,\omega}^{-1}\mathcal{R}_{2}$. Furthermore, since $\widetilde{\mathcal{F}}^{-1}_{2,\omega} \mathcal{R}_{2}$ can be viewed as a low-rank modification of $\mathcal{F}_{2,\omega}^{-1} \mathcal{R}_{2}+\mathcal{Q}_\tau$ via $\mathcal{P}_\tau$ (which has a bounded $\ell_2$-norm), most of the eigenvalues of $\widetilde{\mathcal{F}}^{-1}_{2,\omega} \mathcal{R}_{2}$ are clustered around those of $\mathcal{F}_{2,\omega}^{-1} \mathcal{R}_{2}+\mathcal{Q}_\tau$. In summary, along with Remark \ref{remark41}, when $\omega=\omega^{\star}=\sqrt{\mathcal{O}(\Delta t^{2})+1}$ or $\omega =1$, the eigenvalues of $\widetilde{\mathcal{F}}^{-1}_{2,\omega} \mathcal{R}_{2}$ are nearly clustered within a neighborhood of 1 with a radius $\mathcal{O}(\epsilon + \Delta t)$.
\end{remark}

\section{Numerical experiments}\label{exp}
In this section, we present an extensive array of numerical outcomes pertaining to the solution of 1D and 2D RFNSE utilizing preconditioned GMRES (PGMRES) methods, which incorporate the novel sine-transform-based preconditioner. The discretization of RFNSE relies on a linearly implicit difference scheme that necessitates both the initial value at the onset of time and an approximate value of second or higher order at the subsequent time level. The initial value is directly provided by the initial condition, whereas the latter is derived through the application of a second-order scheme; see \cite{WDL2013JCP} for the 1D case, and \cite{2022dis2D} for the 2D case.

To illustrate the effectiveness and efficiency of the sine-transform-based preconditioner (simply denoted by $\mathcal{F}_{\tau}=\widetilde{\mathcal{F}}_{d,\omega}$ for both $d=1,2$), numerical experiments employing PGMRES methods alongside circulant-based preconditioners are presented for comparative analysis. We denote by $\mathcal{F}_{C}$ the circulant-based preconditioner simply reading that
$$\mathcal{F}_{C}=\frac{1}{2\omega}(\omega I+\mathcal{C}_{d})(\omega I+\mathcal{D}_{d}),\quad d=1,2,$$
where
$\mathcal{C}_{1}=\bigl[\begin{smallmatrix}
	0 & C \\
	-C & 0
\end{smallmatrix}\bigr]$, $\mathcal{C}_{2}=\bigl[\begin{smallmatrix}
	0 & \hat{C} \\
	-\hat{C} & 0
\end{smallmatrix}\bigr]$, $\hat{C}=I\otimes C+C\otimes I$ with $C$ being the Strang's circulant approximation \cite{CRH1989SISSC} of $T_{1}$.

In all the numerical experiments, the linear system at the 2nd time level of the discrete RFNSE is used for testing. The initial guess of the PGMRES method is set as the zero vector. We denote by `$\tau$-GMRES' the PGMRES method with $\mathcal{F}_{\tau}$, `C-GMRES' the PGMRES method with $\mathcal{F}_{C}$, `CPU' the computing time in seconds, and `IT' the number of iterations. Furthermore, the PGMRES method is running without restart, and the stopping criterion is chosen as the $\ell_2$-norm relative residual of the tested linear system reduced below $10^{-8}$ or the number of iterations exceeding $2000$.

\subsection{1D RFNSE with attractive nonlinearities}
We focus on the following truncated 1D RFNSE with an attractive nonlinearity
\begin{align}\label{1DRFSE}
	\imath u_t+\frac{\partial^{\alpha}u(x,t)}{\partial|x|^{\alpha}} + \rho \vert u \vert^2u=0, \quad  -20\le x \le 20,\ 0 < t \le \mathbf{t},
\end{align}
 with the initial and Dirichlet boundary conditions
\begin{align}\nonumber
	u(x,0)=\text{sech}(x) \  e^{2\imath x},\ -20\le x \le 20; \quad u(-20,t)=u(20,t)=0, \ 0 < t \le \mathbf{t}.
\end{align}

\begin{figure}[htbp]
	\centering
	\includegraphics[scale=0.45]{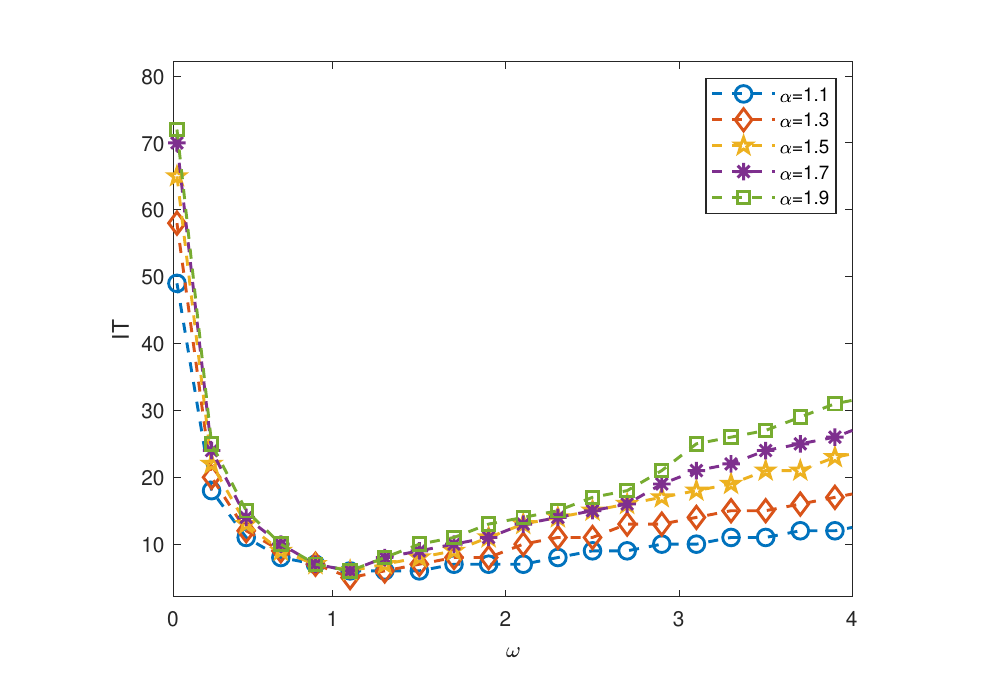}
	\caption{The curves of IT versus the iteration parameter $\omega\in (0,4]$ of $\tau$-GMRES when $\alpha=1.1:0.2:1.9$, $\rho=2$, $M=6400$ and $N=200$.}\label{1Dwit}
\end{figure}

Figure \ref{1Dwit} shows the curves of IT of $\tau$-GMRES versus of the iteration parameter $\omega\in (0,4]$ when $\rho=2$, $M=6400, N=200$ with different fractional orders $\alpha=1.1:0.2:1.9$.
We observe that when $\omega$ approaches to $0$, IT increases rapidly. Meanwhile, as $\omega$ gets larger, IT quickly reaches its minimum and then grows slowly. Thus, the convergence of $\tau$-GMRES is insensitive to the parameter $\omega$ away from 0. In addition, for all tested values of $\alpha$, the optimal parameters are closely to the right of $1$.

\begin{figure}[htbp]
	\centering
	\subfloat{\includegraphics[scale=0.35]{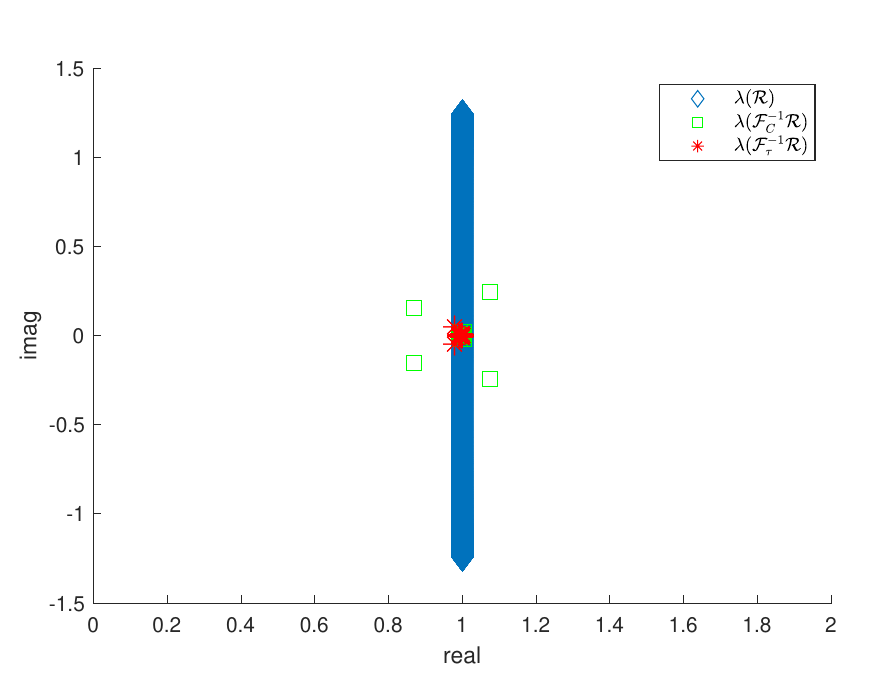}}
	\subfloat{\includegraphics[scale=0.35]{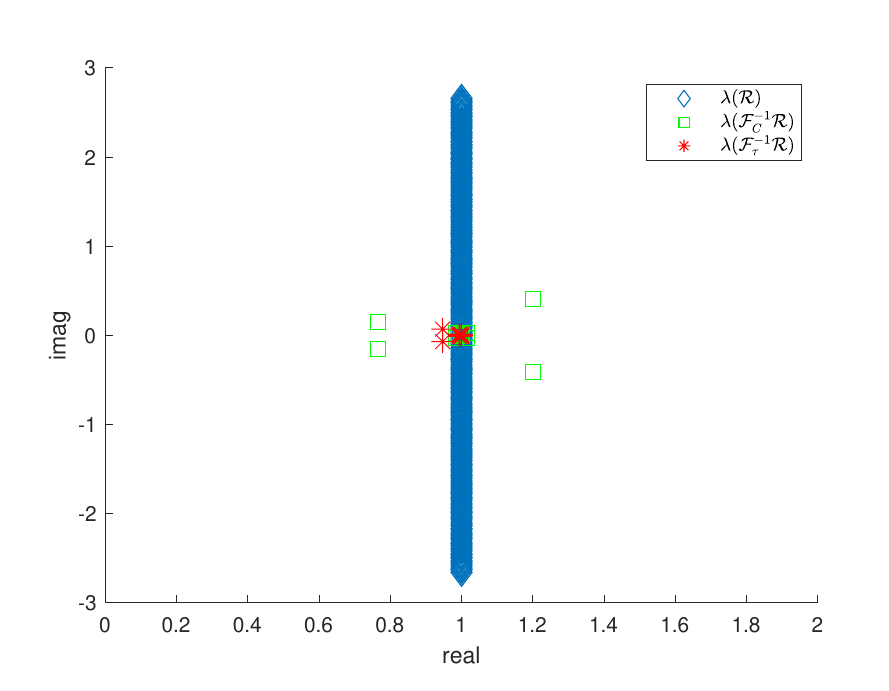}}
	\caption{The eigenvalue distribution of $\mathcal{R}$, ${\mathcal{F}_{\tau}^{-1}\mathcal{R}}$, ${\mathcal{F}_{C}^{-1}\mathcal{R}}$ when $\alpha=1.1$, $\rho=2$, $N=200$, $M=1600$ (left) and $M=3200$ (right).}
	\label{eig11}
\end{figure}

\begin{figure}[htbp]
	\centering
	\subfloat{\includegraphics[scale=0.35]{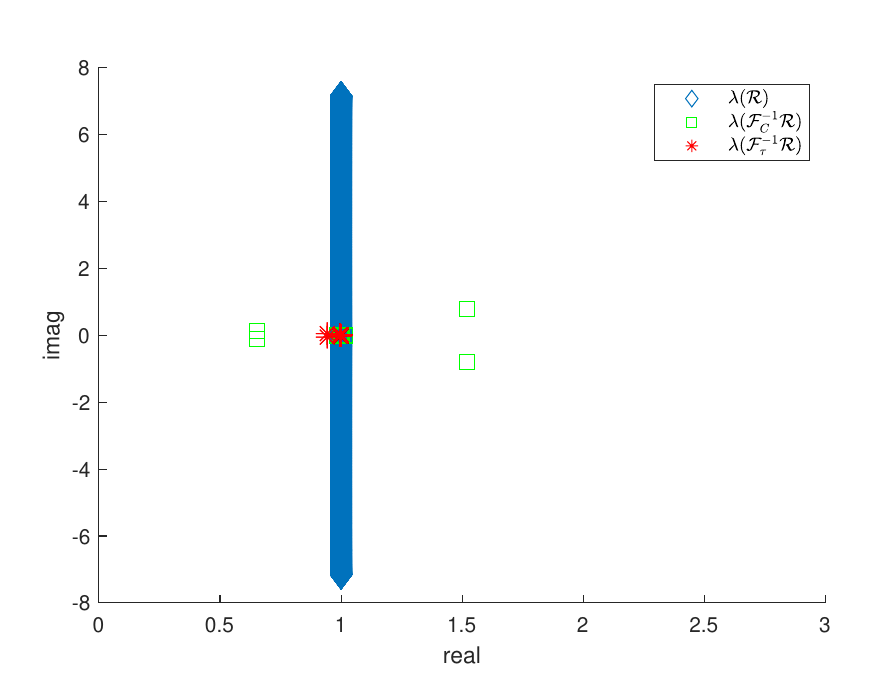}}
	\subfloat{\includegraphics[scale=0.35]{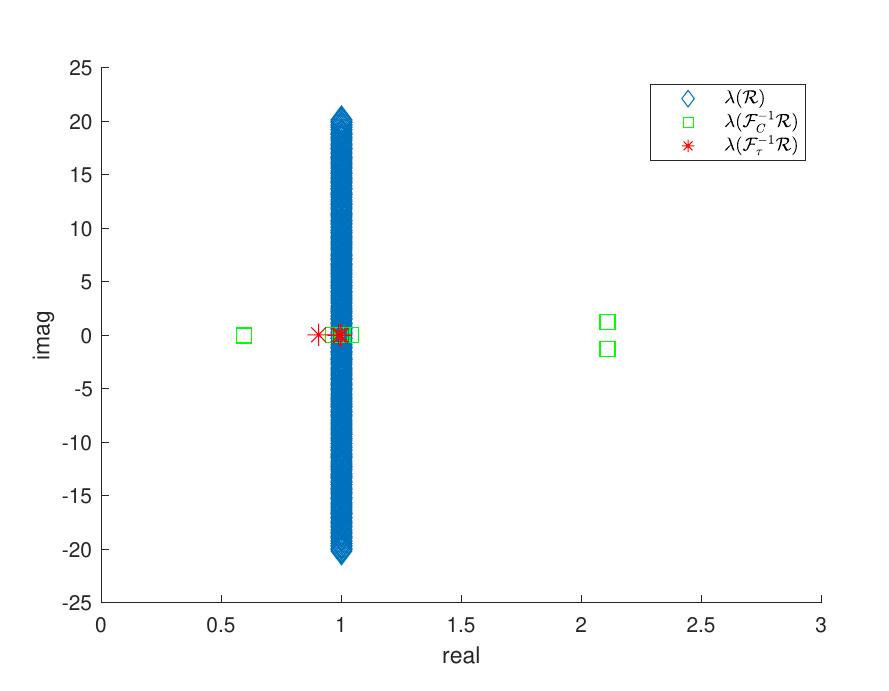}}
	\caption{The eigenvalue distribution of $\mathcal{R}$, ${\mathcal{F}_{\tau}^{-1}\mathcal{R}}$, ${\mathcal{F}_{C}^{-1}\mathcal{R}}$ when $\alpha=1.5$, $\rho=2$, $N=200$, $M=1600$ (left) and $M=3200$ (right).}\label{eig15}
\end{figure}

\begin{figure}[htbp]
	\centering
	\subfloat{\includegraphics[scale=0.35]{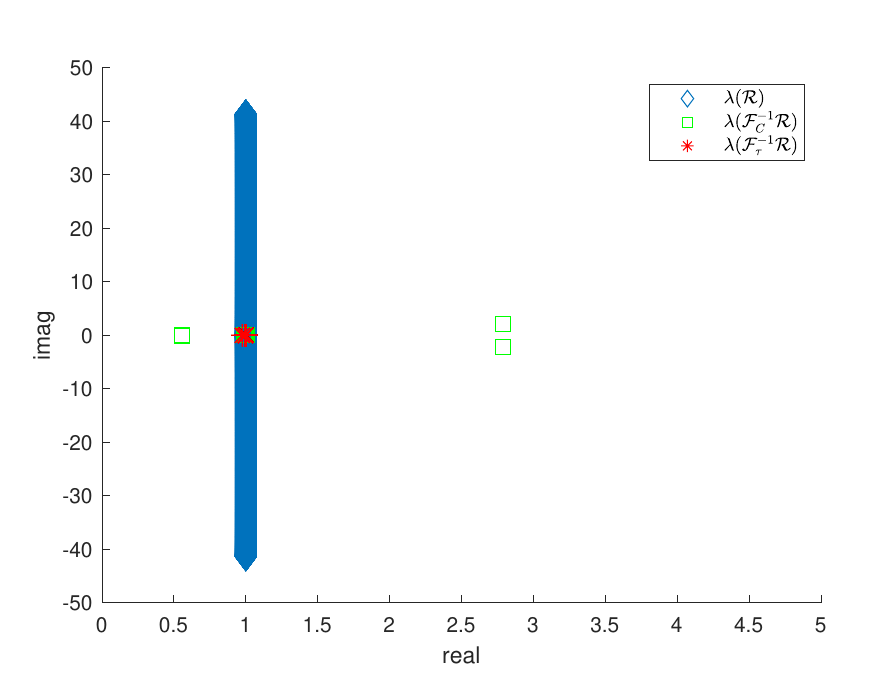}}
	\subfloat{\includegraphics[scale=0.35]{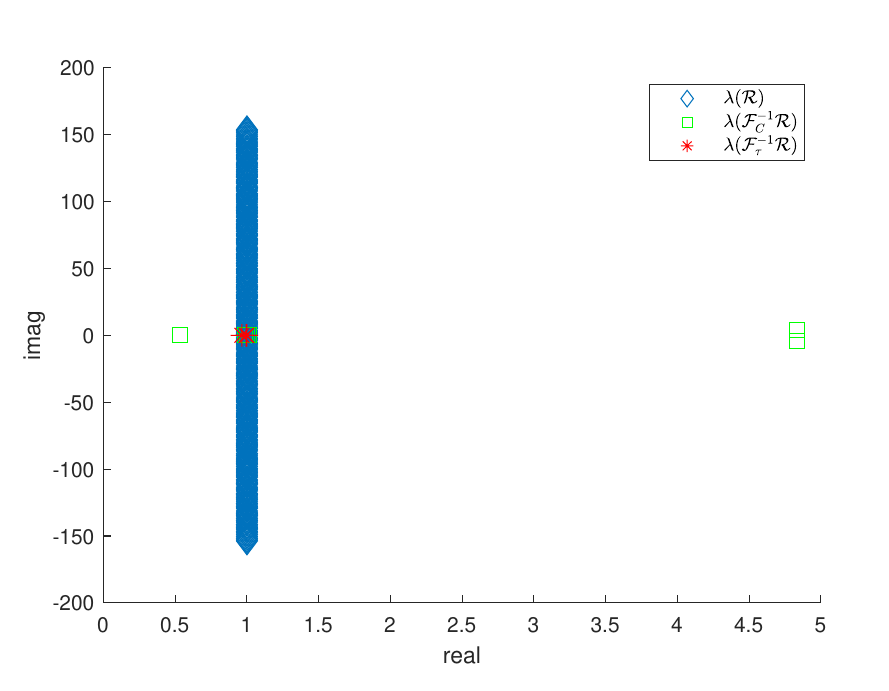}}
	\caption{The eigenvalue distribution of $\mathcal{R}$, ${\mathcal{F}_{\tau}^{-1}\mathcal{R}}$, ${\mathcal{F}_{C}^{-1}\mathcal{R}}$ when $\alpha=1.9$, $\rho=2$, $N=200$, $M=1600$ (left) and $M=3200$ (right).}\label{eig19}
\end{figure}

Figures \ref{eig11}-\ref{eig19} illustrate the eigenvalue distributions of the system matrix $\mathcal{R}$, the circulant-based preconditioned system matrix ${\mathcal{F}_{C}^{-1}\mathcal{R}}$, and the sine-transform-based preconditioned system matrix ${\mathcal{F}_{\tau}^{-1}\mathcal{R}}$. In each figure, the $x$-axis represents the real parts of the eigenvalues, while the $y$-axis represents the imaginary parts. The left side corresponds to $M=1600$, and the right side corresponds to
$M=3200$. In these figures, the real parts of the eigenvalues of $\mathcal{R}$ are consistently 1, while the imaginary parts are widely distributed. For instance, in Figure \ref{eig15}, when $M=3200$, the imaginary part of the eigenvalues of $\mathcal{R}$ range from $-20$ to $20$. It is observed that as $\alpha$ and $M$ increase, the distribution of the imaginary parts of the eigenvalues of $\mathcal{R}$ becomes broader. In contrast, the real parts of the eigenvalues of ${\mathcal{F}_{C}^{-1}\mathcal{R}}$ and ${\mathcal{F}_{\tau}^{-1}\mathcal{R}}$ are clustered around 1, with their imaginary parts being close to 0 (e.g., ranging from $-2$ to $2$ in Figure \ref{eig15}, when $M=3200$). This indicates that the eigenvalues of ${\mathcal{F}_{C}^{-1}\mathcal{R}}$ and ${\mathcal{F}_{\tau}^{-1}\mathcal{R}}$ are more tightly grouped compared to those of $\mathcal{R}$, with
${\mathcal{F}_{\tau}^{-1}\mathcal{R}}$ exhibiting an even tighter clustering than
${\mathcal{F}_{C}^{-1}\mathcal{R}}$. Furthermore, as the spatial mesh size $M$ increases from $1600$ to $3200$, the eigenvalue distributions of ${\mathcal{F}_{\tau}^{-1}\mathcal{R}}$ remain largely unchanged. This observation indirectly suggests that the convergence of $\tau$-GMRES is independent of the spatial mesh size, which is confirmed by Tables \ref{1DITCPU12}-\ref{1DITCPU18}.

The experimental results for IT and CPU of $\tau$-GMRES, C-GMRES and GMRES are presented in Tables \ref{1DITCPU12}-\ref{1DITCPU18}. Here, we set $\alpha=1.2:0.2:1.8$, $\rho=2$, $M=6400, 12800,25600,51200,102400$, and $N=200$. The symbol `-' indicates that the method did not converge within the prescribed maximum iteration count or ran out of memory. As shown in Tables \ref{1DITCPU12}-\ref{1DITCPU18}, GMRES consistently exhibits the highest IT and CPU among all tested methods. Furthermore, as $\alpha$ and $M$ increase, GMRES may fail to converge within the specified maximum iteration count or may ran out of memory. In contrast, both C-GMRES and $\tau$-GMRES demonstrate significantly lower IT and CPU compared to GMRES. Notably, IT of $\tau$-GMRES is smaller than that of C-GMRES, underscoring the superiority of $\tau$-GMRES. Unlike C-GMRES and GMRES, IT of $\tau$-GMRES remains unchanged as the spatial mesh size $M$ increases. This independence on spatial mesh size makes $\tau$-GMRES an excellent tool for solving the linear systems arising from the 1D RFNSE.
\begin{table}[htbp]
	\setlength{\abovecaptionskip}{0pt}
	\setlength{\belowcaptionskip}{10pt} \centering{
		\caption{\label{1DITCPU12} 
			IT and CPU of $\tau$-GMRES, C-GMRES and GMRES when $\rho=2$ and $\alpha=1.2$.}\small
		\begin{tabular}{lcccccc}\specialrule{0em}{2pt}{2pt}\hline\specialrule{0em}{2pt}{2pt}
			 &  $\tau$-GMRES &  & C-GMRES & & GMRES &
			\\\cmidrule(l){2-3}\cmidrule(l){4-5}\cmidrule(l){6-7}
			& IT & CPU & IT & CPU & IT & CPU
			\\\specialrule{0em}{1pt}{1pt}\hline\specialrule{0em}{3pt}{3pt}
			$M=6400$ &  6 & 5.23E-02 & 	8 &	3.47E-02 &	317 & 2.72E+00
			\\\specialrule{0em}{3pt}{3pt}
			$M=12800$ &  6 &	4.57E-02 & 8 &	6.94E-02 & 648 &	2.34E+01
			\\\specialrule{0em}{3pt}{3pt}
			$M=25600$ &  6 & 8.82E-02 & 8 & 1.29E-01 & 1375 & 1.64E+02
			\\\specialrule{0em}{3pt}{3pt}
			$M=51200$ & 6	& 1.70E-01 & 8	& 2.76E-01 & - & -
			\\\specialrule{0em}{3pt}{3pt}
			$M=102400$ & 6 & 3.42E-01 & 8 & 4.60E-01 & - & -
			\\\specialrule{0em}{3pt}{3pt}\hline
	\end{tabular}}
\end{table}

\begin{table}[htbp]
	\setlength{\abovecaptionskip}{0pt}
	\setlength{\belowcaptionskip}{10pt} \centering{
		\caption{\label{1DITCPU14} 
			IT and CPU of $\tau$-GMRES, C-GMRES and GMRES when $\rho=2$ and $\alpha=1.4$.}\small
		\begin{tabular}{lcccccc}\specialrule{0em}{2pt}{2pt}\hline\specialrule{0em}{2pt}{2pt}
			&  $\tau$-GMRES &  & C-GMRES & & GMRES &
			\\\cmidrule(l){2-3}\cmidrule(l){4-5}\cmidrule(l){6-7}
			& IT & CPU & IT & CPU & IT & CPU
			\\\specialrule{0em}{1pt}{1pt}\hline\specialrule{0em}{3pt}{3pt}
			$M=6400$ &  6 & 2.81E-02 & 	8 &	2.75E-02 &	1299 & 4.57E+01
			\\\specialrule{0em}{3pt}{3pt}
			$M=12800$ &  6 &	4.01E-02 & 8 &	7.06E-02 & - &	-
			\\\specialrule{0em}{3pt}{3pt}
			$M=25600$ &  6 & 8.32E-02 & 9 & 1.41E-01 & - & -
			\\\specialrule{0em}{3pt}{3pt}
			$M=51200$ & 6	& 1.64E-01 & 9	& 3.29E-01 & - & -
			\\\specialrule{0em}{3pt}{3pt}
			$M=102400$ & 6 & 3.16E-01 & 8 & 5.16E-01 & - & -
			\\\specialrule{0em}{3pt}{3pt}\hline
	\end{tabular}}
\end{table}

\begin{table}[htbp]
	\setlength{\abovecaptionskip}{0pt}
	\setlength{\belowcaptionskip}{10pt} \centering{
		\caption{\label{1DITCPU16} 
			IT and CPU of $\tau$-GMRES, C-GMRES and GMRES when $\rho=2$ and $\alpha=1.6$.}\small
		\begin{tabular}{lcccccc}\specialrule{0em}{2pt}{2pt}\hline\specialrule{0em}{2pt}{2pt}
			&  $\tau$-GMRES &  & C-GMRES & & GMRES &
			\\\cmidrule(l){2-3}\cmidrule(l){4-5}\cmidrule(l){6-7}
			& IT & CPU & IT & CPU & IT & CPU
			\\\specialrule{0em}{1pt}{1pt}\hline\specialrule{0em}{3pt}{3pt}
			$M=6400$ &  6 & 1.94E-02 & 	9 &	2.80E-02 &	- & -
			\\\specialrule{0em}{3pt}{3pt}
			$M=12800$ &  6 &	4.00E-02 & 9 &	7.81E-02 & - &	-
			\\\specialrule{0em}{3pt}{3pt}
			$M=25600$ &  6 & 8.74E-02 & 10 & 1.56E-01 & - & -
			\\\specialrule{0em}{3pt}{3pt}
			$M=51200$ & 6	& 1.59E-01 & 11	& 4.11E-01 & - & -
			\\\specialrule{0em}{3pt}{3pt}
			$M=102400$ & 6 & 3.30E-01 & 10 & 5.79E-01 & - & -
			\\\specialrule{0em}{3pt}{3pt}\hline
	\end{tabular}}
\end{table}

\begin{table}[htbp]
	\setlength{\abovecaptionskip}{0pt}
	\setlength{\belowcaptionskip}{10pt} \centering{
		\caption{\label{1DITCPU18} 
			IT and CPU of $\tau$-GMRES, C-GMRES and GMRES when $\rho=2$ and $\alpha=1.8$.}\small
		\begin{tabular}{lcccccc}\specialrule{0em}{2pt}{2pt}\hline\specialrule{0em}{2pt}{2pt}
			&  $\tau$-GMRES &  & C-GMRES & & GMRES &
			\\\cmidrule(l){2-3}\cmidrule(l){4-5}\cmidrule(l){6-7}
			& IT & CPU & IT & CPU & IT & CPU
			\\\specialrule{0em}{1pt}{1pt}\hline\specialrule{0em}{3pt}{3pt}
			$M=6400$ &  6 & 1.80E-02 & 	11 &	3.53E-02 &	- & -
			\\\specialrule{0em}{3pt}{3pt}
			$M=12800$ &  6 &	3.93E-02 & 12 &	1.06E-01 & - &	-
			\\\specialrule{0em}{3pt}{3pt}
			$M=25600$ &  6 & 8.27E-02 & 14 & 2.20E-01 & - & -
			\\\specialrule{0em}{3pt}{3pt}
			$M=51200$ & 6	& 1.68E-01 & 16	& 5.85E-01 & - & -
			\\\specialrule{0em}{3pt}{3pt}
			$M=102400$ & 6 & 3.24E-01 & 14 & 8.29E-01 & - & -
			\\\specialrule{0em}{3pt}{3pt}\hline
	\end{tabular}}
\end{table}

\begin{figure}[htbp]
	\centering
	\subfloat{\includegraphics[scale=0.35]{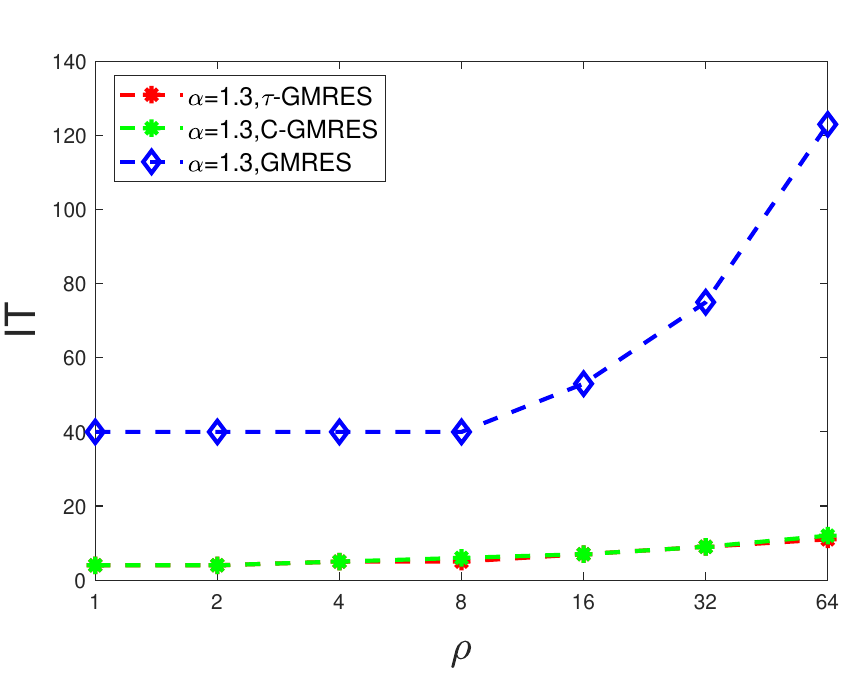}}
	\subfloat{\includegraphics[scale=0.35]{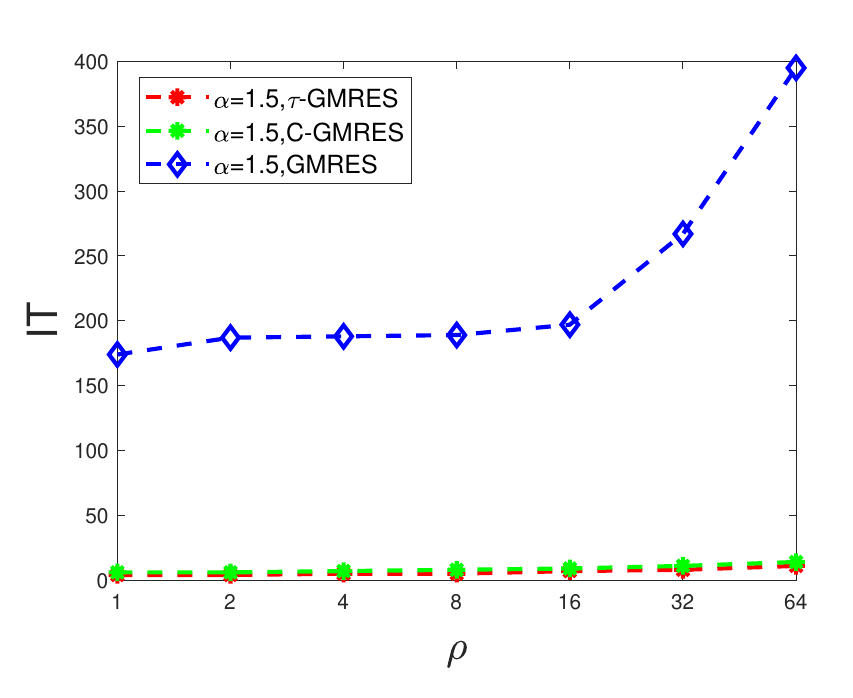}}
	\subfloat{\includegraphics[scale=0.35]{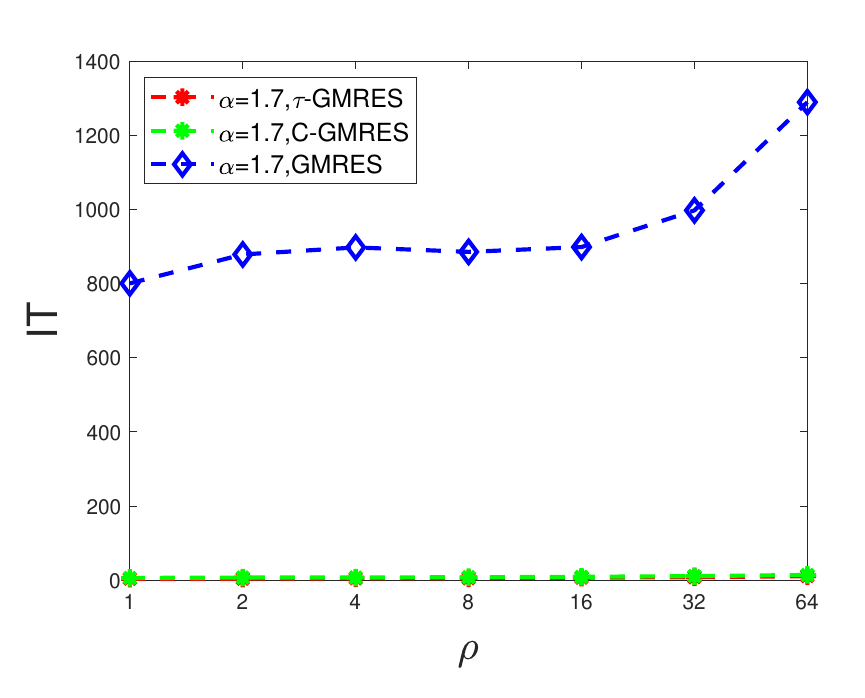}}
	\caption{The curves of IT of $\tau$-GMRES, C-GMRES and GMRES versus the nonlinear term parameter $\rho$ when $\alpha=1.3:0.2:1.7$, $M=6400$, $N=200$.}\label{1DrhoIT}
\end{figure}

Figure \ref{1DrhoIT} illustrates the impact of the strength of the nonlinear term (controlled by the parameter $\rho$) on IT of GMRES, C-GMRES, and $\tau$-GMRES, with $\alpha=1.3:0.2:1.7$, $M=6400$, and $N=200$. As the parameter $\rho$ increases from $1$ to $64$, IT of GMRES grows rapidly. In contrast, IT of C-GMRES and $\tau$-GMRES rises slowly. In summary, the preconditioned GMRES methods with the preconditioners $\mathcal{F}_{C}$ and $\mathcal{F}_{\tau}$ demonstrate resilience to the strength of the nonlinear term, highlighting the reliability of the new preconditioned GMRES method.

\begin{figure}[htbp]
	\centering
	\includegraphics[scale=0.35]{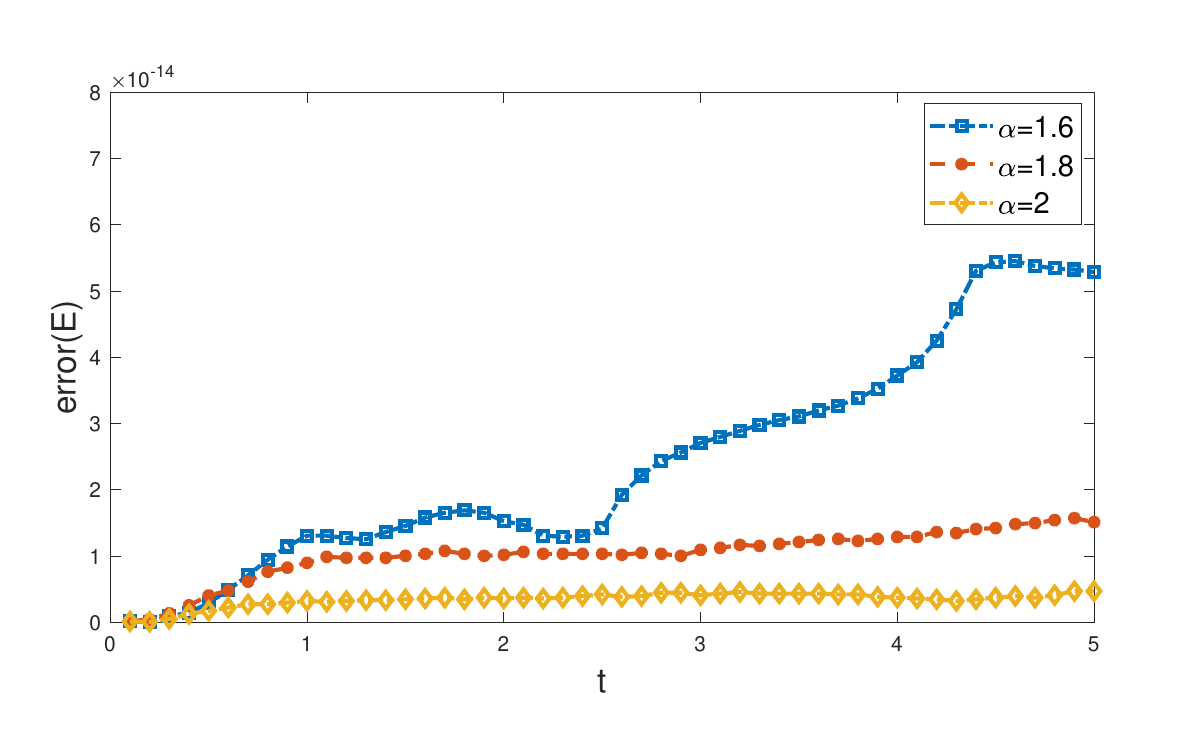}
	\caption{The relative errors of the discrete energy 
with $\rho=2$, $h=0.2$ and $\Delta t=0.001$.}\label{1DEerror}
\end{figure}

Table \ref{1DQE} and Figure \ref{1DEerror} present the relative errors of the discrete mass and energy for $\rho=2$ and various fractional orders $\alpha$. In Table \ref{1DQE}, We set the spatial step size $h=0.2$ and the temporal step size $\Delta t=0.05$, while in Figure \ref{1DEerror}, $h=0.2$ and $\Delta t=0.001$. At each time level, we utilize the $\tau$-GMRES method to solve the linear system (\ref{discretizedCNLSMaxForm}), with the stopping criterion defined as the $\ell_2$-norm relative residual of the tested system being reduced below $10^{-15}$. As shown in Table \ref{1DQE} and Figure \ref{1DEerror}, the relative errors of the discrete mass and energy are very small, indicating that the $\tau$-GMRES method effectively preserves the conservation properties of the LICD scheme \cite{WDL2014JCP}.

\begin{table}[htbp]
	\setlength{\abovecaptionskip}{0pt}
	\setlength{\belowcaptionskip}{10pt} \centering{\caption{\label{1DQE}
			The relative errors of the discrete mass 
with $\rho=2$, $h=0.2$ and $\Delta t=0.05$.}\small
		\begin{tabular}{lllllllllllllllllll}\specialrule{0em}{2pt}{2pt}\hline\specialrule{0em}{2pt}{2pt}
			& & & & &  $t=1$ & & & & $t=2$ & & & & $t=3$ & & & & $t=4$ &
			\\\specialrule{0em}{1pt}{1pt}\hline\specialrule{0em}{3pt}{3pt}
			& $\alpha=1.4$ & & & & 2.2216e-16 & & & & 2.2216e-16  & & & & 2.2216e-16 & & & & 5.5540e-16 &
			\\\specialrule{0em}{3pt}{3pt}
			& $\alpha=1.7$ & & & & 1.1110e-16 & & & &  1.1110e-16  & & & &  5.5548e-16 & & & &  5.5548e-16 &
			\\\specialrule{0em}{3pt}{3pt}
			& $\alpha=1.9$ & & & & 0 & & & & 4.4444e-16 & & & &  4.4444e-16 & & & &  0 &
			\\\specialrule{0em}{3pt}{3pt}
			&  $\alpha=2$ & & & &  3.3335e-16 & & & &  0  & & & &  3.3335e-16  & & & &  0 &
			\\\specialrule{0em}{3pt}{3pt}\hline
	\end{tabular}}
\end{table}

\subsection{2D RFNSE with attractive nonlinearities}
We consider the truncated 2D RFNSE with an attractive nonlinearity as follows
\begin{align}\label{2DRFSE}
	\imath u_t+\frac{\partial^{\alpha}u(\bbb{x},t)}{\partial|\bbb{x}|^{\alpha}} + \rho\vert u \vert^2u=0,\quad \bbb{x}=(x,y)\in\Omega,\ t\in(0,\mathbf{t}],
\end{align}
with the initial and Dirichlet boundary conditions
\begin{align}\nonumber
	u(x,y,0)=\frac{2}{\sqrt{\pi}}\text{exp}[-(x^{2}+y^{2})],\ (x,y)\in\Omega;\quad u(x,y,t)=0,\ (x,y)\in\partial\Omega.
\end{align}
We take the truncated domain as $\Omega=[\mbox{-5},\mbox{5}]^{2}$.

\begin{figure}[htbp]
	\centering
	\includegraphics[scale=0.40]{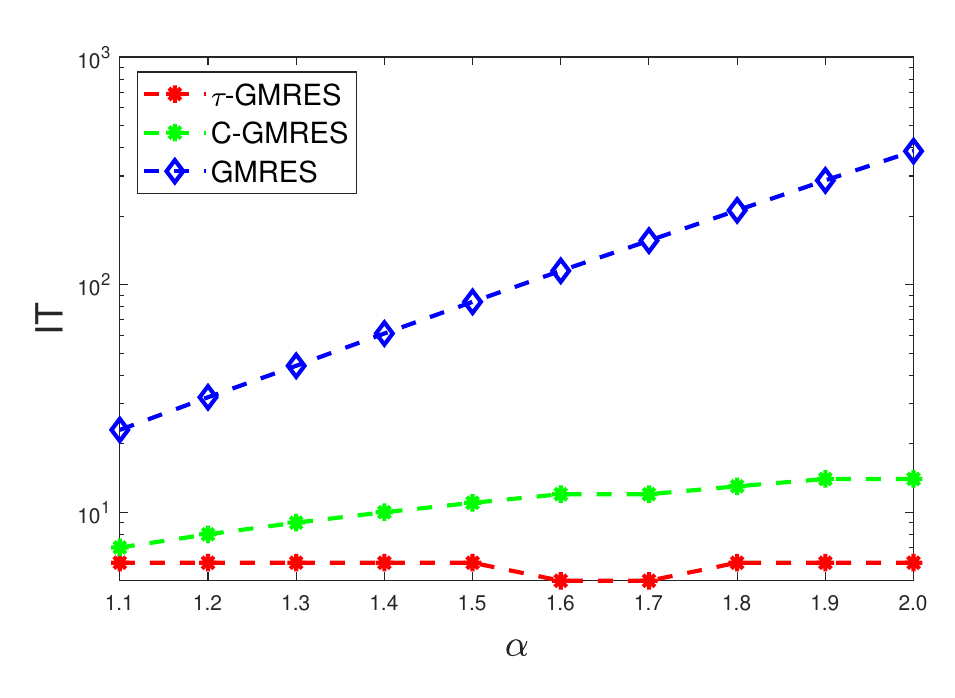}
	\caption{IT of $\tau$-GMRES, C-GMRES and GMRES as a function of the fractional order $\alpha=1.1:0.1:2.0$ when $\rho=1$, $h=1/8$, $\Delta t=1/20$.}\label{2Dalpit}
\end{figure}

Figure \ref{2Dalpit} illustrates IT of $\tau$-GMRES, C-GMRES and GMRES as a function of the fractional order $\alpha$. We set the spatial step size to $h=1/8$, the temporal step size to $\Delta t=1/20$, the nonlinear term parameter $\rho=1$, and the fractional order to $\alpha=1.1:0.1:2.0$. The red, green, and blue lines represent the IT variation curves of $\tau$-GMRES, C-GMRES, and GMRES, respectively. As shown in Figure \ref{2Dalpit}, as the fractional order $\alpha$ increases, IT of GMRES rises very rapidly. In contrast, IT of C-GMRES increases linearly with a very small slope, while IT of $\tau$-GMRES remains almost constant. This indicates that both $\tau$-GMRES and C-GMRES are not sensitive to the fractional order $\alpha$, with $\tau$-GMRES demonstrating superior performance over C-GMRES in terms of iteration count.

Tables \ref{2DITCPU12}-\ref{2DITCPU18} present the numerical results regarding IT and CPU of $\tau$-GMRES,C-GMRES and GMRES when $\rho=1$, $\alpha=1.2:0.2:1.8$ and $h=1/32$, $1/64$, $1/128$, $1/256$, $1/512$. It is important to note that the symbol `$n$' represents the scale of the complex linear system (\ref{discretizedRFSE2DMaxForm}). As shown in Tables \ref{2DITCPU12}-\ref{2DITCPU18}, GMRES exhibit the highest IT and CPU among all tested methods. When either the fractional order $\alpha$ or the matrix-size scale $n$ increases, GMRES may fail to converge within the the allowed maximum iteration count or run out of memory. This highlights the challenges GMRES faces in large-scale scenarios. In contrast, IT and CPU of $\tau$-GMRES and C-GMRES are significantly lower than those of GMRES, with $\tau$-GMRES outperforming C-GMRES in both metrics. Furthermore, as the fractional order $\alpha$ and matrix-size scale $n$ increase, IT of C-GMRES grows linearly, while IT of $\tau$-GMRES remains nearly constant. This indicates that $\tau$-GMRES significantly enhances the computational efficiency when solving the complex linear system (\ref{discretizedRFSE2DMaxForm}) and demonstrates independence from the spatial mesh size and the fractional order.

\begin{table}[htbp]
	\setlength{\abovecaptionskip}{0pt}
	\setlength{\belowcaptionskip}{10pt} \centering{
		\caption{\label{2DITCPU12} 
			IT and CPU of $\tau$-GMRES, C-GMRES and GMRES when $\rho=1$ and $\alpha=1.2$.}\small
		\begin{tabular}{lccccccc}\specialrule{0em}{2pt}{2pt}\hline\specialrule{0em}{2pt}{2pt}
			& &  $\tau$-GMRES &  & C-GMRES & & GMRES &
			\\\cmidrule(l){3-4}\cmidrule(l){5-6}\cmidrule(l){7-8}
			& $n$& IT & CPU & IT & CPU & IT & CPU
			\\\specialrule{0em}{1pt}{1pt}\hline\specialrule{0em}{3pt}{3pt}
			$h=1/32$ & 102400 & 6 & 5.86E-02 & 	12 &	9.68E-01 &	179 & 2.51E+01
			\\\specialrule{0em}{3pt}{3pt}
			$h=1/64$ & 409600 & 6 &	3.32E+00 & 11 &	5.85E-01 & 351 &	3.41E+02
			\\\specialrule{0em}{3pt}{3pt}
			$h=1/128$ & 1638400& 5 & 1.12E+01 & 11 & 2.62E+01 & - & -
			\\\specialrule{0em}{3pt}{3pt}
			$h=1/256$ & 6553600& 6	& 6.98E+01 & 11	& 1.16E+02 & - & -
			\\\specialrule{0em}{3pt}{3pt}
			$h=1/512$ &26214400 & 6 & 5.02E+02 & - & - & - & -
		\\\specialrule{0em}{3pt}{3pt}\hline
		\end{tabular}}
	\end{table}

\begin{table}[htbp]
	\setlength{\abovecaptionskip}{0pt}
	\setlength{\belowcaptionskip}{10pt} \centering{
		\caption{\label{2DITCPU14} 
			IT and CPU of $\tau$-GMRES, C-GMRES and GMRES when $\rho=1$ and $\alpha=1.4$.}\small
		\begin{tabular}{lccccccc}\specialrule{0em}{2pt}{2pt}\hline\specialrule{0em}{2pt}{2pt}
			& &  $\tau$-GMRES &  & C-GMRES & & GMRES &
			\\\cmidrule(l){3-4}\cmidrule(l){5-6}\cmidrule(l){7-8}
			& $n$& IT & CPU & IT & CPU & IT & CPU
			\\\specialrule{0em}{1pt}{1pt}\hline\specialrule{0em}{3pt}{3pt}
			$h=1/32$ & 102400 & 6 &	5.00E-01 & 	14 &	1.18E+00 &	443 &	1.18E+02
			\\\specialrule{0em}{3pt}{3pt}
			$h=1/64$ & 409600 & 6 &	3.30E+00 & 15 &	9.48E+00 & 1069 &	2.54E+03
			\\\specialrule{0em}{3pt}{3pt}
			$h=1/128$ & 1638400& 5 & 1.20E+01 & 16 & 3.96E+01 & - & -
			\\\specialrule{0em}{3pt}{3pt}
			$h=1/256$ & 6553600& 5	& 5.85E+01 & 16	& 1.88E+02 & - & -
			\\\specialrule{0em}{3pt}{3pt}
			$h=1/512$ &26214400 & 6 & 4.98E+02 & - & - & - & -
			\\\specialrule{0em}{3pt}{3pt}\hline
	\end{tabular}}
\end{table}

\begin{table}[htbp]
	\setlength{\abovecaptionskip}{0pt}
	\setlength{\belowcaptionskip}{10pt} \centering{
		\caption{\label{2DITCPU16} 
			IT and CPU of $\tau$-GMRES, C-GMRES and GMRES when $\rho=1$ and $\alpha=1.6$.}\small
		\begin{tabular}{lccccccc}\specialrule{0em}{2pt}{2pt}\hline\specialrule{0em}{2pt}{2pt}
			& &  $\tau$-GMRES &  & C-GMRES & & GMRES &
			\\\cmidrule(l){3-4}\cmidrule(l){5-6}\cmidrule(l){7-8}
			& $n$& IT & CPU & IT & CPU & IT & CPU
			\\\specialrule{0em}{1pt}{1pt}\hline\specialrule{0em}{3pt}{3pt}
			$h=1/32$ & 102400 & 6 &	5.34E-01 & 	19 &	1.58E+00 &	1085 &	6.23E+02
			\\\specialrule{0em}{3pt}{3pt}
			$h=1/64$ & 409600 & 6 &	4.08E+00 & 22 &	1.45E+01 & - &	-
			\\\specialrule{0em}{3pt}{3pt}
			$h=1/128$ & 1638400& 5 & 1.22E+01 & 26 & 6.57E+01 & - & -
			\\\specialrule{0em}{3pt}{3pt}
			$h=1/256$ & 6553600& 5	& 5.96E+01 & 32	& 4.16E+02 & - & -
			\\\specialrule{0em}{3pt}{3pt}
			$h=1/512$ &26214400 & 6 & 4.87E+02 & - & - & - & -
			\\\specialrule{0em}{3pt}{3pt}\hline
	\end{tabular}}
\end{table}

\begin{table}[htbp]
	\setlength{\abovecaptionskip}{0pt}
	\setlength{\belowcaptionskip}{10pt} \centering{
		\caption{\label{2DITCPU18} 
			IT and CPU of $\tau$-GMRES, C-GMRES and GMRES when $\rho=1$ and $\alpha=1.8$.}\small
		\begin{tabular}{lccccccc}\specialrule{0em}{2pt}{2pt}\hline\specialrule{0em}{2pt}{2pt}
			& &  $\tau$-GMRES &  & C-GMRES & & GMRES &
			\\\cmidrule(l){3-4}\cmidrule(l){5-6}\cmidrule(l){7-8}
			& $n$& IT & CPU & IT & CPU & IT & CPU
			\\\specialrule{0em}{1pt}{1pt}\hline\specialrule{0em}{3pt}{3pt}
			$h=1/32$ & 102400 & 6 &	4.94E-01 & 	26 &	2.22E+00 &	- &	-
			\\\specialrule{0em}{3pt}{3pt}
			$h=1/64$ & 409600 & 6 &	4.10E+00 & 34 &	2.26E+01 & - &	-
			\\\specialrule{0em}{3pt}{3pt}
			$h=1/128$ & 1638400& 6 & 1.65E+01 & 47 & 1.31E+02 & - & -
			\\\specialrule{0em}{3pt}{3pt}
			$h=1/256$ & 6553600& 6	& 7.38E+01 & 95	& 2.53E+03 & - & -
			\\\specialrule{0em}{3pt}{3pt}
			$h=1/512$ &26214400 & 6 & 5.14E+02 & - & - & - & -
			\\\specialrule{0em}{3pt}{3pt}\hline
	\end{tabular}}
\end{table}

\begin{figure}[htbp]
	\centering
	\subfloat{\includegraphics[scale=0.40]{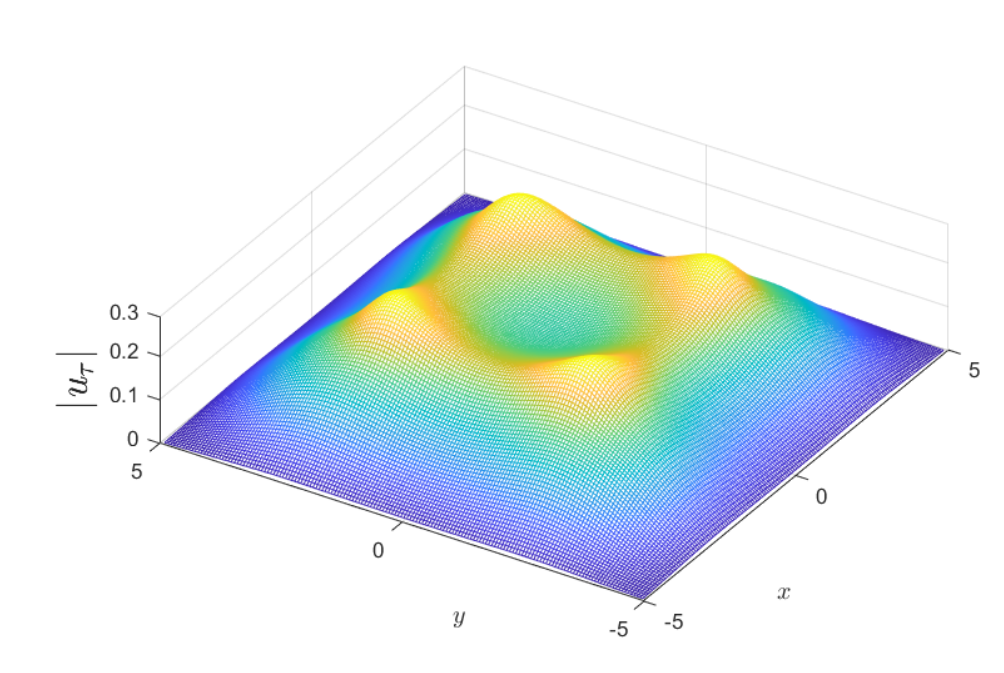}}
	\subfloat{\includegraphics[scale=0.40]{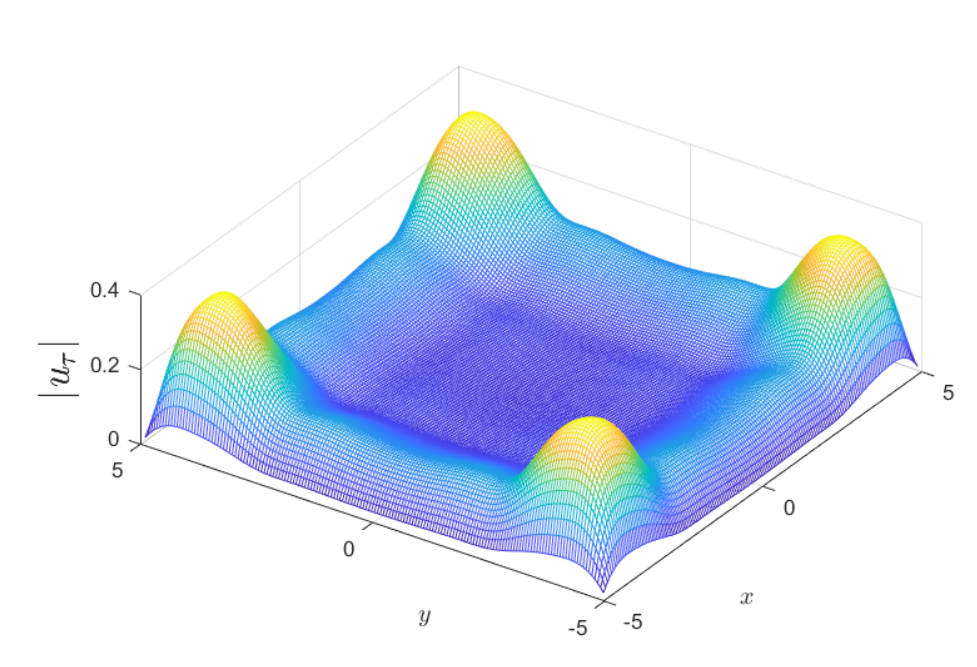}}
	\caption{The profile of the numerical solution $u_{\tau}$ for $\rho=1$, $\alpha=1.3$ and $h=\Delta t=1/16$ at $\mathbf{t}=2$ (left), $\mathbf{t}=4$ (right).}\label{2Dsol13}
\end{figure}

\begin{figure}[htbp]
	\centering
	\subfloat{\includegraphics[scale=0.40]{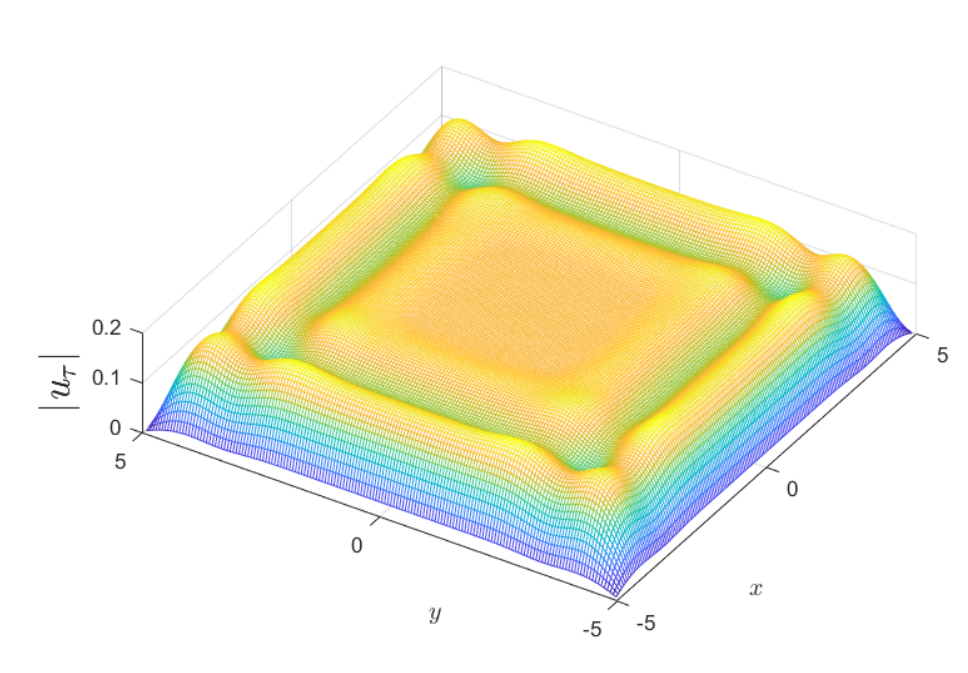}}
	\subfloat{\includegraphics[scale=0.40]{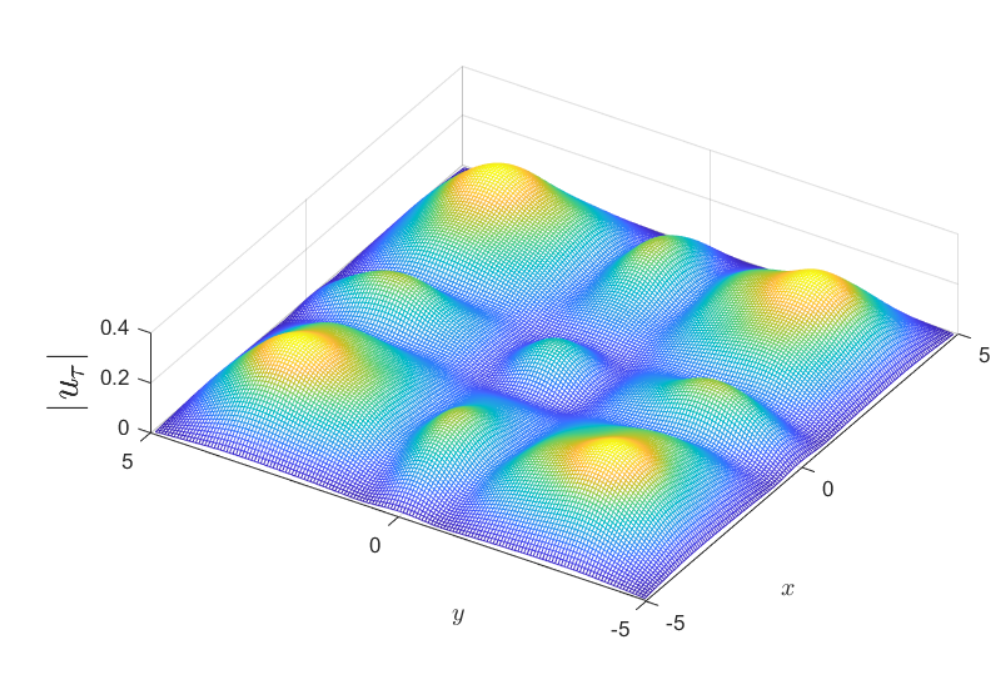}}
	\caption{The profile of the numerical solution $u_{\tau}$ for $\rho=1$, $\alpha=1.7$ and $h=\Delta t=1/16$ at $\mathbf{t}=2$ (left), $\mathbf{t}=4$ (right).}\label{2Dsol17}
\end{figure}

\begin{figure}[htbp]
	\centering
	\subfloat{\includegraphics[scale=0.40]{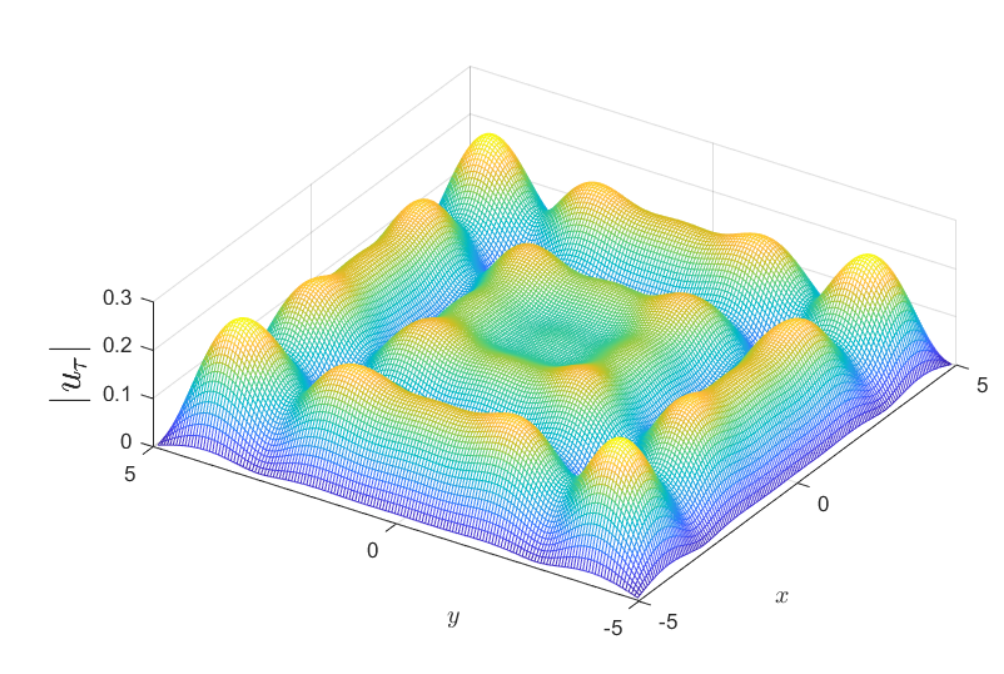}}
	\subfloat{\includegraphics[scale=0.40]{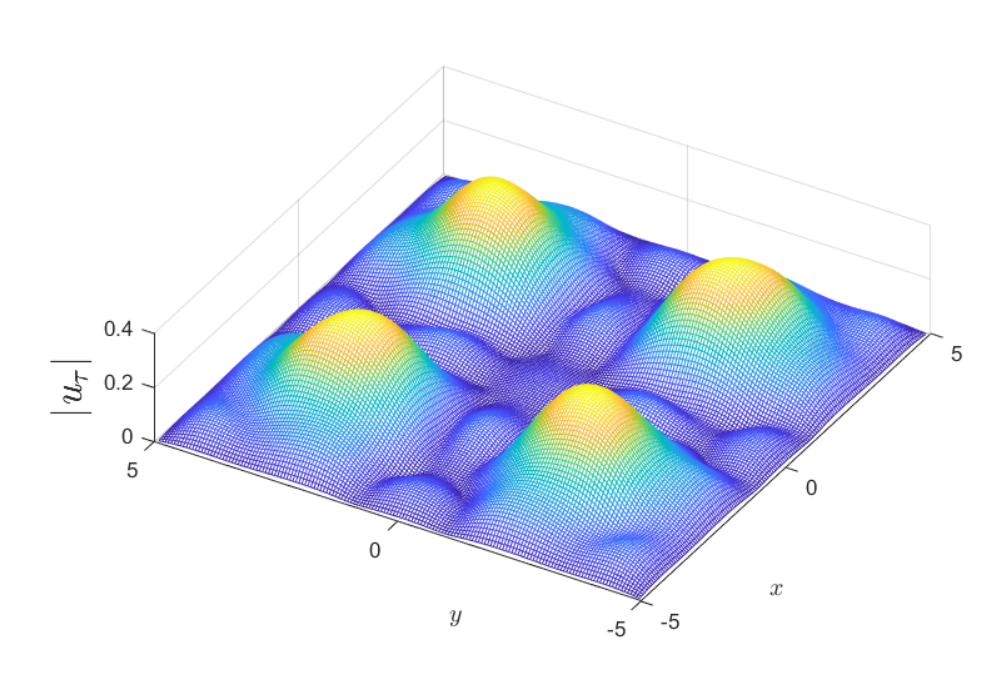}}
	\caption{The profile of the numerical solution $u_{\tau}$ for $\rho=1$, $\alpha=1.95$ and $h=\Delta t=1/16$ at $\mathbf{t}=2$ (left), $\mathbf{t}=4$ (right).}\label{2Dsol195}
\end{figure}

Figures \ref{2Dsol13}-\ref{2Dsol195} illustrate the profiles of the numerical solution $u_{\tau}$ obtained using $\tau$-GMRES when $\rho=1$ and $h=\Delta t=1/16$. The values of the fractional order $\alpha$ are referenced from \cite{2022dis2D}. In each figure, the left side corresponds to the time $\mathbf{t}=2$, while the right side corresponds to $\mathbf{t}=4$. The results indicate that the value of $\alpha$ has a significant impact on the shape of the numerical solution $u_{\tau}$. Specifically, as $\alpha$ increases, the numerical solution exhibits greater heterogeneity, a phenomenon that becomes particularly pronounced after the system has evolved over a longer duration.

\begin{figure}[htbp]
	\centering
	\subfloat{\includegraphics[scale=0.40]{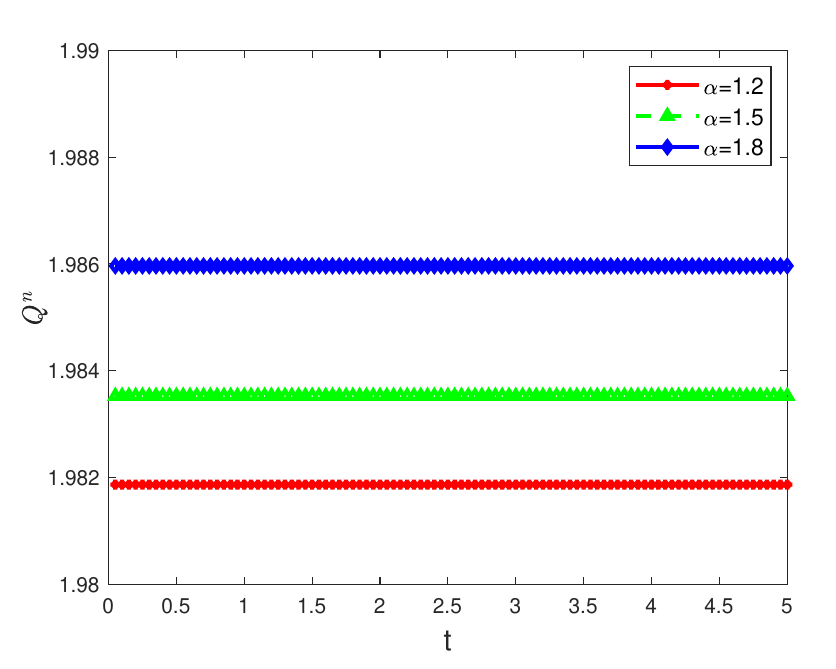}}
	\subfloat{\includegraphics[scale=0.40]{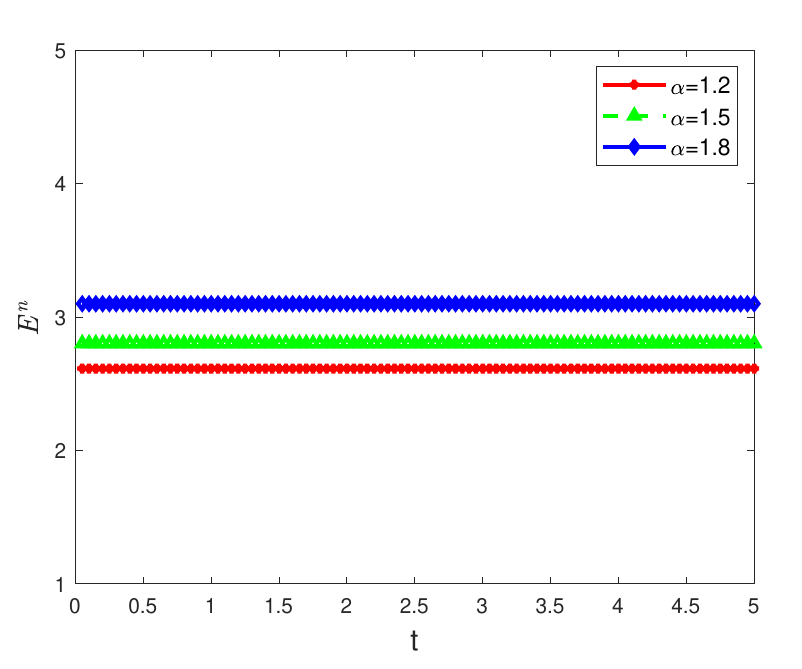}}
	\caption{The evolution of the discrete mass $Q^n$ (left) and energy $E^n$ (right) when $h=\Delta t=1/20$, $\rho=1$, and $\alpha=1.2:0.3:1.8$. }\label{2DQE}
\end{figure}

Figure \ref{2DQE} illustrates the evolution of discrete mass $Q^n$ and energy $E^n$ defined in \cite{2022dis2D}. We set $h=\Delta t=1/20$, the final time $\mathbf{t}=5$, the nonlinear term parameter $\rho=1$, and the fractional order $\alpha=1.2:0.3:1.8$. The results in Figure \ref{2DQE} demonstrate that $\tau$-GMRES effectively preserves the conservation properties of the TLID scheme \cite{2022dis2D}. Furthermore, while the value of discrete mass remains relatively weak correlation with changes in the fractional order $\alpha$, the discrete energy shows a relatively strong correlation with $\alpha$.

\section{Concluding remarks}\label{clu}

This paper presents a TBAN iteration method for solving indefinite complex linear systems derived from the discretization of Schrödinger equations with Riesz fractional derivatives and attractive nonlinear terms. These linear systems exhibit complex symmetry, indefiniteness, and a $d$-level Toeplitz-plus-diagonal structure. Theoretical analysis shows that the TBAN iteration method possesses unconditional convergence and the parameter-free property, with the optimal parameter approximately equal to 1. By combining the aforementioned iteration method with the sine-transform-based preconditioning technique, we construct a novel preconditioner to accelerate the convergence speed of the GMRES method. The eigenvalues of the system matrix preconditioned by the new preconditioner demonstrate good clustering, and the convergence behavior of the corresponding PGMRES method is independent of the spatial mesh size and the fractional order, also exhibiting the parameter-free property (i.e., the optimal parameter can be fixed at 1). Regarding future work, variable coefficient problems and non-uniform spatial discretization schemes may result in complex linear systems that do not possess an explicit $d$-level Toeplitz-plus-diagonal structure. Consequently, investigating potential implicit data-sparse structures and integrating other methods with the framework of our approach to develop new preconditioners will be both a valuable and challenging undertaking.

\vspace{2em}
\section*{Acknowledgments}
This work was funded by the National Natural Science Foundation (No. 11101213 and No. 12071215), China.

\end{document}